\documentclass[11pt]{article}

\usepackage{mathptmx}       
\usepackage{amsfonts}
\usepackage{amsbsy}
\SetSymbolFont{letters}{bold}{OML}{cmr}{bx}{n}
\SetSymbolFont{letters}{normal}{OML}{cmr}{m}{n}

\usepackage{helvet}
\usepackage{courier}
\usepackage{type1cm}        

\usepackage{graphicx}
\usepackage{makeidx}
\usepackage{multicol}
\usepackage[bottom]{footmisc}

\usepackage{booktabs}
\usepackage{multirow}
\usepackage[caption=false]{subfig}

\begin{document}
\title{Dynamic Bayesian Combination of Multiple Imperfect Classifiers}
\author{Edwin Simpson$^1$, Stephen Roberts$^1$, Ioannis Psorakis$^1$, Arfon Smith$^2$\\ \\
1. Department of Engineering Science, University of Oxford, UK.\\
Emails: \{edwin, sjrob, psorakis\}@robots.ox.ac.uk\\
2. Zooniverse. Email: arfon@zooniverse.org}
\maketitle

\abstract{Classifier combination methods need to make best use of the outputs of multiple, imperfect classifiers to enable higher accuracy classifications. In many situations, such as when human decisions need to be combined, the base decisions can vary enormously in reliability. A Bayesian approach to such uncertain combination allows us to infer the differences in performance between individuals and to incorporate any available prior knowledge about their abilities when training data is sparse. In this paper we explore Bayesian classifier combination, using the computationally efficient framework of variational Bayesian inference. We apply the approach to real data from a large citizen science project, Galaxy Zoo Supernovae, and show that our method far outperforms other established approaches to imperfect decision combination. We go on to analyse the putative community structure of the decision makers, based on their inferred decision making strategies, and show that natural groupings are formed. Finally 
we present a dynamic Bayesian classifier combination approach and investigate the changes in base classifier performance over time.}

\section{Introduction}

In many real-world scenarios we are faced with the need to aggregate information from cohorts of imperfect decision making agents (\emph{base classifiers}), be they computational or human. Particularly in the case of human agents, we rarely have available to us an indication of how decisions were arrived at or a realistic measure of agent confidence in the various decisions. Fusing multiple sources of information in the presence of uncertainty is optimally achieved using Bayesian inference, which elegantly provides a principled mathematical framework for such knowledge aggregation. In this paper we provide a Bayesian framework for imperfect decision combination, where the base classifications we receive are greedy preferences (i.e. labels with no indication of confidence or uncertainty). The classifier combination method we develop aggregates the decisions of multiple agents, improving overall performance. We present a principled framework in which the use of weak decision makers can be mitigated and in 
which multiple agents, with very different observations, knowledge or training sets, can be combined to provide complementary information. 

The preliminary application we focus on in this paper is a distributed \emph{citizen science} project, in which human agents carry out classification tasks, in this case identifying transient objects from images as corresponding to potential supernovae or not. This application, \emph{Galaxy Zoo Supernovae} \cite{smith_galaxy_2010}, is part of the highly successful \emph{Zooniverse} family of citizen science projects. In this application the ability of our base classifiers can be very varied and there is no guarantee over any individual's performance, as each user can have radically different levels of domain experience and have different background knowledge. As individual users are not overloaded with decision requests by the system, we often have little performance data for individual users (base classifiers). The methodology we advocate provides a scalable, computationally efficient, Bayesian approach to learning base classifier performance thus enabling optimal decision combinations. The approach is 
robust in the presence of uncertainties at all levels and naturally handles missing observations, i.e. in cases where agents do not provide any base classifications. We develop extensions to allow for dynamic, sequential inference, through which we track information regarding the base classifiers. Through the application of social network analysis we can also observe behavioural patterns in the cohort of base classifiers.

The rest of this paper is organised as follows. In the remainder of the Introduction we briefly describe related work. In Section \ref{sec:IBCCmodel} we present a probabilistic model for independent Bayesian classifier combination, IBCC. Section \ref{sec:VB} introduces the approximate inference method, variational Bayes, and details its application to IBCC. Section \ref{sec:galaxyzooresults} shows an example application for classifier combination, Galaxy Zoo Supernovae, and compares results using different classifier combination methods, including IBCC. In Sections \ref{sec:picomms} and \ref{sec:taskcomms} we investigate how communities of decision makers with similar characteristics can be found using data inferred from Bayesian classifier combination. Section \ref{sec:dynibcc} presents an extension to independent Bayesian classifier combination that models the changing performance of individual decision makers. Using this extension, Section \ref{sec:ContributorDynamics} examines the dynamics of individuals 
from our example application, while Sections \ref{sec:picommsDynamics} and \ref{sec:taskcommsDynamics} show how communities of decision makers change over time. Finally, Section \ref{sec:discussion} discusses future directions for this work.

\subsection{Related Work}

Previous work has often focused on aggregating expert decisions in fields such as medical diagnosis \cite{dawid_maximum_1979}. In contrast, \emph{crowdsourcing} uses novice human agents to perform tasks that would be too difficult or expensive to process computationally or using experts \cite{bloodgood_using_2010,quinn_crowdflow:_2010}. The underlying problem of fusing labels from multiple classifications has been dealt with in various ways and a review of the common methods is given by \cite{tulyakov_review_2008}. The choice of method typically depends on the type of labels we can obtain from agents (e.g. binary, continuous), whether we can manipulate agent performance, and whether we can also access input features. Weighted majority and weighted sum algorithms are popular methods that account for differing reliability in the base classifiers; an analysis of their performance is given by \cite{littlestone_weighted_2002}. Bayesian model combination \cite{monteith_turning_2011} provides a theoretical basis 
for soft-selecting from a space of combination functions. In most cases it outperforms Bayesian model averaging, which relies on one base classifier matching the data generating model. A well-founded model that learns the combination function directly was defined by \cite{ghahramani_bayesian_2003}, giving a Bayesian treatment to a model first presented in \cite{dawid_maximum_1979}. A similar model was also investigated by \cite{raykar_learning_2010} with extensions to learn a classifier from expert labels rather than known ground truth labels. Both papers assume that base classifiers have constant performance, a problem that we address later in this paper.

\section{Independent Bayesian Classifier Combination} \label{sec:IBCCmodel}

Here we present a variant of Independent Bayesian Classifier Combination (IBCC), originally defined in \cite{ghahramani_bayesian_2003}. The model assumes conditional independence between base classifiers, but performed as well as more computationally intense dependency modelling methods also given by \cite{ghahramani_bayesian_2003}. 

We are given a set of data points indexed from $1$ to $N$, where the $i$th data point has a true label $t_{i}$ that we wish to infer. We assume $t_i$ is generated from a multinomial distribution with the probabilities of each class denoted by $\boldsymbol{\kappa:}p(t_{i}=j|\boldsymbol{\kappa})=\kappa_{j}$. True labels may take values $t_i = {1...J}$, where $J$ is the number of true classes. We assume there are $K$ base classifiers, which produce a set of discrete outputs $\vec{c}$ with values $l={1..L}$, where $L$ is the number of possible outputs. The output $c^{(k)}_i$ from classifier $k$ for data point $i$ is assumed to be generated from a
 multinomial distribution dependent on the true label, with parameters \mbox{ $\boldsymbol{\pi}^{(k)}_j:p(c_{i}^{(k)}=l|t_{i}=j,\boldsymbol{\pi}_j^{(k)})=\pi_{jl}^{(k)}$ }. This model places minimal requirements on the type of classifier output, which need not be probabilistic and could be selected from an arbitrary number of discrete values, indicating, for example, greedy preference over a set of class labels. Parameters $\boldsymbol{\pi}_j^{(k)}$ and $\boldsymbol{\kappa}$ have Dirichlet prior distributions with hyperparameters $\boldsymbol{\alpha}_{0,j}^{(k)}=[\alpha_{0,j1}^{(k)},...,\alpha_{0,jL}^{(k)}]$ and $\boldsymbol{\nu}=[\nu_{01},...\nu_{0J}]$ respectively. We refer to the set of $\boldsymbol{\pi}_j^{(k)}$ for all base classifiers and all classes as $ \boldsymbol{\Pi}=\left\{ \boldsymbol{\pi}^{(k)}_j | j=1...J, k=1...K\right\} $. Similarly, for the hyperparameters we use $ \mathbf{A}_0=\left\{ \boldsymbol{\alpha_0}^{(k)}_j | j=1...J, k=1...K\right\} $

\begin{figure}[ht]
\centering{}\includegraphics[clip=false,width=0.6\textwidth,bb=0 0 500 300]{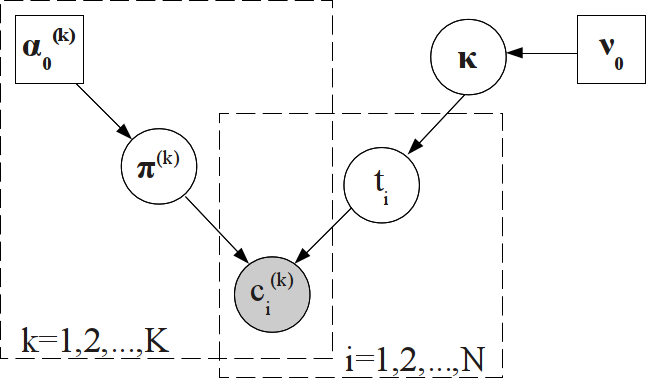}
\caption{Graphical Model for IBCC. The shaded node represents observed values, circular
nodes are variables with a distribution and square nodes are variables instantiated
with point values.
\label{fig:graphical_model}}
\end{figure}

The joint distribution over all variables for the IBCC model is 
\begin{equation}
p(\boldsymbol{\kappa},\boldsymbol{\Pi},\vec{t},\vec{c}|\mathbf{A}_0,\boldsymbol{\nu})=\prod_{i=1}^{N}\{\kappa_{t_{i}}\prod_{k=1}^{K}\pi^{(k)}_{t_{i},c_{i}^{(k)}}\}p(\boldsymbol{\kappa}|\boldsymbol{\nu})p(\boldsymbol{\Pi}|\mathbf{A})\!,\label{eq:IBCC}
\end{equation}
The graphical model for IBCC is shown in Figure \ref{fig:graphical_model}.
A key feature of IBCC is that $\boldsymbol{\pi}^{(k)}$ represents a \emph{confusion matrix} that quantifies the decision-making abilities of \emph{each} base classifier. This potentially allows us to ignore or retrain poorer classifiers and assign expert decision makers to data points that are highly uncertain. Such efficient selection of base classifiers is vitally important when there is a cost to obtaining an output from a base classifier, for example, a financial payment to a decision maker per decision, or when the processing bandwidth of base classifiers is limited. The confusion matrices $\boldsymbol{\Pi}$ in IBCC also allow us to predict any missing classifier outputs in $\vec{c}$, so that we can naturally handle cases where only \emph{partially observed} agents make decisions.

The IBCC model assumes independence between the rows in $\boldsymbol{\pi}^{(k)}$, i.e. the probability of each classifier's outputs is dependent on the \emph{true
label class}. In some cases it may be reasonable to assume that performance over one label class may be correlated with performance in another; indeed methods such as \emph{weighted majority} \cite{littlestone_weighted_2002} make this tacit assumption. However, we would argue that this is not universally the case, and IBCC makes no such strong assumptions.

The model here represents a simplification of that proposed in \cite{ghahramani_bayesian_2003}, which places exponential hyperpriors over $\boldsymbol{A_{0}}$. The exponential distributions are not conjugate to the Dirichlet distributions over $\boldsymbol{\Pi}$,
so inference using \emph{Gibbs Sampling} \cite{geman_stochastic_1984} requires an expensive adaptive rejection sampling step \cite{gilks_adaptive_1992} for $\mathbf{A}$ and the variational Bayesian solution is intractable. The conjugate prior to the Dirichlet is non-standard and its normalisation constant is not in closed form \cite{lefkimmiatis_bayesian_2009}, so cannot be used. We therefore alter the model, to use point values for $\boldsymbol{A_{0}}$, as in other similar models \cite{choudrey03,bishopbooko,penny_dynamic_2002}. The hyperparameter values of $\boldsymbol{A_{0}}$ can hence be chosen to represent any prior level of uncertainty in the values of the agent-by-agent confusion matrices, $\boldsymbol{\Pi}$, and can be regarded as pseudo-counts of prior observations, offering a natural method to include any prior knowledge and a methodology to extend the method to sequential, on-line environments.

\section{Variational Bayesian IBCC} \label{sec:VB}

The goal of the combination model is to perform inference for the unknown variables $\vec{t},$ $\boldsymbol{\Pi},$ and $\boldsymbol{\kappa}$. The inference technique proposed in \cite{raykar_learning_2010} is maximum a posteriori (MAP) estimation, while \cite{ghahramani_bayesian_2003} suggests a full Bayesian treatment using Gibbs Sampling \cite{geman_stochastic_1984}. While the latter provides some theoretical guarantee of accuracy given the proposed model, it is often very slow to converge and convergence is difficult to ascertain. In this paper we consider the use of principled approximate Bayesian methods, namely \emph{variational Bayes} (VB) \cite{attias_advances_2000} as this allows us to replace non-analytic marginal integrals in the original model with analytic updates in the \emph{sufficient statistics} of the variational approximation. This produces a model that iterates rapidly to a solution in a computational framework which can be seen as a Bayesian generalisation of the \emph{Expectation-
Maximization} (EM) algorithm \cite{dempster_maximum_1977}.

\subsection{Variational Bayes}

Given a set of observed data $\boldsymbol{X}$ and a set of latent variables and parameters $\boldsymbol{Z}$, the goal of variational Bayes (VB) is to find a tractable approximation $q(\boldsymbol{Z})$ to the posterior distribution $p(\boldsymbol{Z}|\boldsymbol{X})$ by minimising the KL-divergence \cite{kullback_information_1951} between the approximate distribution
and the true distribution \cite{attias_advances_2000, fox_tutorial_2011}. We can write the log of the model evidence $p(\boldsymbol{X})$ as 
\begin{eqnarray}
\ln p(\mathbf{X}) & = & \int q(\mathbf{Z})\ln{\frac{p(\mathbf{X,Z})}{q(\mathbf{Z})}}\mathrm{d}\mathbf{Z}-\int q(\mathbf{Z})\ln{\frac{p(\mathbf{Z}|\mathbf{X})}{q(\mathbf{Z})}}\mathrm{d}\mathbf{Z}\\ \nonumber
 & = & L(q)-\mathrm{KL}(q||p). \label{eq:VBEQ}
\end{eqnarray}
As $q(\mathbf{Z})$ approaches $p(\mathbf{Z}|\mathbf{X})$, the KL-divergence disappears and the lower bound $L(q)$ on $\ln p(\mathbf{X})$ is maximised. Variational Bayes selects a restricted form of $q(\mathbf{Z})$ that is tractable to work with, then seeks the distribution within this restricted form that minimises the KL-divergence. A common restriction is to partition $\mathbf{Z}$ into groups of variables, then assume $q(\mathbf{Z})$ factorises into functions of single groups:
\begin{equation}
  q(\mathbf{Z})=\prod_{i=1}^{M}q_{i}(\mathbf{Z}_{i}) .
\end{equation}
For each factor $q_{i}(\mathbf{Z}_{i})$ we then seek the optimal solution $q_{i}^{*}(\mathbf{Z}_{i})$ that minimises the KL-divergence. Consider partitions of variables $\mathbf{Z}_i$ and $\mathbf{\bar{Z}}_i$, where $\mathbf{\bar{Z}}_i=\{\mathbf{Z}_j|j\neq i,j=1...M\}$. \emph{Mean field theory} \cite{parisi1988statistical} shows that we can find an optimal factor $q_{i}^{*}(\mathbf{Z}_{i})$ from the conditional distribution $p(\mathbf{Z}_{i}|\mathbf{X},\mathbf{\bar{Z}})$ by taking the expectation over all the other factors ${j|j\neq i,j=1...M}$. We can therefore write the log of the optimal factor $\ln q_{i}^{*}(\mathbf{Z}_{i})$ as the expectation with respect to all other factors of the log of the joint distribution over all variables plus a normalisation constant:
\begin{equation}
\ln q_{i}^{*}(\mathbf{Z}_{i})=\mathbb{E}_{q(\mathbf{\bar{Z}})}[\ln p(\mathbf{X},\mathbf{Z})]+\mathrm{const}.\label{eq:generalOptimalFactor}
\end{equation}
In our notation, we take the expectation with respect to the variables in the subscript. In this case, $\mathbb{E}_{q(\mathbf{\bar{Z}})}[...]$ indicates that we take an expectation with respect to all factors except $q(\mathbf{Z}_i)$. This expectation is implicitly conditioned on the observed data, $\mathbf{X}$, which we omit from the notation for brevity.

We can evaluate these optimal factors iteratively by first initialising all factors, then updating each in turn using the expectations with respect to the current values of the other factors. Unlike Gibbs sampling, each iteration is guaranteed to increase the lower bound on the log-likelihood, $L(q)$, converging to a (local) maximum in a similar fashion to standard EM algorithms. If the factors $q_{i}^{*}(\mathbf{Z}_{i})$ are exponential family distributions \cite{bishopbooko}, as is the case for the IBCC method we present in the next section, the lower bound is convex with respect to each factor $q_{i}^{*}(\mathbf{Z}_{i})$ and $L(q)$ will converge to a \emph{global}
maximum of our approximate, factorised distribution. In practice, once the optimal factors $q_{i}^{*}(\mathbf{Z}_{i})$ have converged to within a given tolerance, we can approximate the distribution of the unknown variables and calculate their expected values.

\subsection{Variational Equations for IBCC} \label{sec:VBIBCC}

To provide a variational Bayesian treatment of IBCC, VB-IBCC, we first propose the form for our variational distribution, $q(\mathbf{Z})$, that factorises between the parameters and latent variables.
\begin{eqnarray}
q(\boldsymbol{\kappa},\vec{t},\boldsymbol{\Pi}) & = & q(\vec{t})q(\boldsymbol{\kappa},\boldsymbol{\Pi})\label{eq:IBCC-1}
\end{eqnarray}

This is the only assumption we must make to perform VB on this model; the forms of the factors arise from our model of IBCC. We can use the joint distribution in Equation (\ref{eq:IBCC}) to find the optimal factors $q^{*}(\vec{t})$ and $q^{*}(\boldsymbol{\kappa},\boldsymbol{\Pi})$ in the form given by Equation (\ref{eq:generalOptimalFactor}). For the true labels we have
\begin{eqnarray}
\ln q^{*}(\vec{t}) & = & \mathbb{E}_{\boldsymbol{\kappa,\Pi}}[\ln p(\boldsymbol{\kappa},\vec{t},\boldsymbol{\Pi},\vec{c})]+\mathrm{const}.
\end{eqnarray}

We can rewrite this into factors corresponding to independent data points, with any terms not involving $t_{i}$ being absorbed into the normalisation constant. To do this we define $\rho_{ij}$ as 
\begin{eqnarray}
\ln\rho_{ij} & = &  \mathbb{E}_{\kappa_j,\boldsymbol{\pi_j}}[\ln p(\kappa_j,t_i,\boldsymbol{\pi_j},\vec{c})] \\ \nonumber
 & = & \mathbb{E}_{\boldsymbol{\kappa}}[\ln\kappa_{j}]+\sum_{k=1}^{K}\mathbb{E}_{\boldsymbol{\pi}_j}[\ln\pi_{j,c_{i}^{(k)}}^{(k)}]\label{eq:VB_lnqstar_ti}
\end{eqnarray}
then we can estimate the probability of a true label, which also gives its expected value:
\begin{equation}
q^{*}(t_{i} = j) =  \mathbb{E}_{\vec{t}}[t_{i} = j] = \frac{\rho_{ij}}{\sum_{\iota=1}^J \rho_{i\iota} } \;.
\end{equation}
To simplify the optimal factors in subsequent equations, we define expectations with respect to $\vec{t}$ of the number of occurrences of each true class, given by
\begin{equation}
N_{j}=\sum_{i=1}^{N}\mathbb{E}_{\vec{t}}[t_{i}=j] \;,
\label{eq:classCount}
\end{equation}
and the counts of each classifier decision $c_{i}^{(k)}=l$ given the true label $t_{i}=j$, by
\begin{equation}
N_{jl}^{(k)}=\sum_{i=1}^{N}\delta_{c_{i}^{(k)}l}\mathbb{E}_{\vec{t}}[t_{i}=j]
\label{eq:decisionCount}
\end{equation}
where $\delta_{c_{i}^{(k)}l}$ is unity if $c_{i}^{(k)}=l$ and zero otherwise.

For the parameters of the model we have the optimal factors given by:
\begin{eqnarray}
\ln q^{*}(\boldsymbol{\kappa,\Pi}) & = & \mathbb{E}_{\vec{t}}[\ln p(\boldsymbol{\kappa,t,\Pi,}\vec{c})]+\mathrm{const}\label{eq:IBCC-3-2}\\
 & = & \mathbb{E}_{\vec{t}}\left[\sum_{i=1}^{N}\ln {\kappa}_{t_{i}}+\sum_{i=1}^{N}\sum_{k=1}^{K}\ln\pi_{t_{i},c_{i}^{(k)}}^{(k)}\right]+\ln p(\boldsymbol{\kappa}|\boldsymbol{\nu}_0)\\ \nonumber
 &  & +\ln p(\boldsymbol{\Pi}|\mathbf{A}_0)+\mathrm{const}.
\end{eqnarray}
In Equation (\ref{eq:IBCC-3-2}) terms involving $\boldsymbol{\kappa}$ and terms involving each confusion matrix in $\boldsymbol{\Pi}$ are separate, so we can factorise $q^{*}(\boldsymbol{\kappa,\pi})$ further into
\begin{equation}
q^{*}(\boldsymbol{\kappa,}\boldsymbol{\Pi})=q^{*}(\boldsymbol{\kappa})\prod_{k=1}^{K}\prod_{j=1}^{J}q^{*}\left(\boldsymbol{\pi}_{j}^{(k)}\right).\label{eq:IBCC-3-2-1}
\end{equation}
In the IBCC model (Section \ref{sec:IBCCmodel}) we assumed a Dirichlet prior for $\boldsymbol{\kappa}$, which gives us the optimal factor
\begin{eqnarray}
\ln q^{*}(\boldsymbol{\kappa}) & = & \mathbb{E}_{\vec{t}}\left[\sum_{i=1}^{N}\ln\kappa_{t_{i}}\right]+\ln p(\boldsymbol{\kappa}|\boldsymbol{\nu}_0)+\mathrm{const}\label{eq:IBCC-3-2-1-1}\\
 & = & \sum_{j=1}^{J}N_{j}\ln\kappa_{j}+\sum_{j=1}^{J}\left(\nu_{0,j}-1\right)\ln\kappa_{j}+\mathrm{const}.
\end{eqnarray}
Taking the exponential of both sides, we obtain a posterior Dirichlet distribution of the form
\begin{equation}
q^{*}(\boldsymbol{\kappa})\propto\mathrm{Dir}(\boldsymbol{\kappa}|\nu_{1},...,\nu_{J})\label{eq:qStar_P}
\end{equation}
where $\boldsymbol{\nu}$ is updated in the standard manner by adding the data counts to the prior counts $\nu_{0}$:
\begin{equation}
\nu_{j}=\nu_{0,j}+N_{j}.
\end{equation}
The expectation of $\ln\boldsymbol{\kappa}$ required to update Equation (\ref{eq:VB_lnqstar_ti}) is therefore:
\begin{equation}
\mathbb{E}\left[\ln\kappa_{j}\right]=\Psi(\nu_{j})-\Psi\left(\sum_{\iota=1}^{J}\nu_{\iota}\right)\label{eq:eLnPj}
\end{equation}
where $\Psi$ is the standard digamma function \cite{abramowitz_handbook_1965}.

For the confusion matrices $\boldsymbol{\pi}_{j}^{(k)}$ the priors are also Dirichlet distributions giving us the factor 
\begin{eqnarray}
\ln q^{*}\left(\boldsymbol{\pi}_{j}^{(k)}\right) & = & \sum_{i=1}^{N}\mathbb{E}_{t_{i}}[t_{i}=j]\ln\pi_{j,c_{i}^{(k)}}^{(k)}+\ln p\left(\boldsymbol{\pi}^{(k)}_j|\boldsymbol{\alpha}^{(k)}_{0,j}\right)+\mathrm{const}\\
 & = & \sum_{l=1}^{L}N_{jl}^{(k)}\ln\pi_{jl}^{(k)}+\sum_{l=1}^{L}\left(\alpha_{0,jl}^{(k)}-1\right)\ln\pi_{jl}^{(k)}+\mathrm{const}.
\end{eqnarray}
Again, taking the exponential gives a posterior Dirichlet distribution of the form
\begin{eqnarray}
q^{*}\left(\boldsymbol{\pi}_{j}^{(k)}\right) & = & \mathrm{Dir}\left(\boldsymbol{\pi}_{j}^{(k)}|\alpha_{j1}^{(k)},...,\alpha_{jL}^{(k)}\right) \label{eq:VB_qstart_pi}
\end{eqnarray}
where $\boldsymbol{\alpha}_{j}^{(k)}$ is updated by adding data counts to prior counts $\alpha_{0,j}^{(k)}$:
\begin{equation}
\alpha_{jl}^{(k)}=\alpha_{0,jl}^{(k)}+N_{jl}^{(k)}. \label{eq:VB_static_alpha_updates}
\end{equation}
The expectation required for Equation (\ref{eq:VB_lnqstar_ti}) is given by
\begin{equation}
\mathbb{E}\left[\ln\pi_{jl}^{(k)}\right] = \Psi\left(\alpha_{jl}^{(k)}\right)-\Psi\left(\sum_{m=1}^{L}\alpha_{jm}^{(k)}\right).\label{eq:expectation_lnPi}
\end{equation}

To apply the VB algorithm to IBCC, we first choose an initial value for all variables $\mathbb{E}[\ln\pi_{jl}^{(k)}]$ and $\mathbb{E}[\ln\kappa_{j}]$ either randomly or by taking the expectations of the variables over their prior distributions (if we have enough domain knowledge to set informative priors). We then iterate over a two-stage procedure similar to the \emph{Expectation-Maximization} (EM) algorithm. In the variational equivalent of the \emph{E-step} we use the current expected parameters, $\mathbb{E}[\ln\pi_{jl}^{(k)}]$ and $\mathbb{E}[\ln\kappa_{j}]$, to update the variational distribution in Equation (\ref{eq:IBCC-1}). First we evaluate Equation (\ref{eq:VB_lnqstar_ti}), then use the result to update the counts $N_{j}$ and $N_{jl}^{(k)}$ according to Equations (\ref{eq:classCount}) and (\ref{eq:decisionCount}). In the variational \emph{M-step}, we update $\mathbb{E}[\ln\pi_{jl}^{(k)}]$ and $\mathbb{E}[\ln\kappa_{j}]$ using Equations (\ref{eq:eLnPj}) and (\ref{eq:expectation_lnPi}).

\subsection{Variational Lower Bound} \label{sec:staticLowerBound}

To check for convergence we can also calculate the lower bound $L(q)$ (see Equation (\ref{eq:VBEQ})), which should always increase after a pair of \emph{E-step} and \emph{M-step} updates. While we could alternatively detect convergence of the expectation over the latent variables, the variational lower bound is a useful sanity check for our derivation of VB-IBCC and for its implementations.

\begin{eqnarray}
L(q) & = &  \int\!\!\! \int\!\!\! \int   q(\vec{t},\boldsymbol{\Pi},\boldsymbol{\kappa})   \ln\frac{p(\vec{c},\vec{t},\boldsymbol{\Pi},\boldsymbol{\kappa}|\mathbf{A}_{0},\boldsymbol{\nu}_{0})}   {q(\vec{t},\boldsymbol{\Pi},\boldsymbol{\kappa})}  \mathrm{d}\vec{t}  \mathrm{d}\boldsymbol{\Pi}  \mathrm{d}\boldsymbol{\kappa} \nonumber \\
 & = & \mathbb{E}_{\vec{t},\boldsymbol{\Pi},\boldsymbol{\kappa}}[\ln p(\vec{c},\vec{t},\boldsymbol{\Pi},\boldsymbol{\kappa}|\mathbf{A}_{0},\boldsymbol{\nu}_{0})]-\mathbb{E}_{\vec{t},\boldsymbol{\Pi},\boldsymbol{\kappa}}[\ln q(\vec{t},\boldsymbol{\Pi},\boldsymbol{\kappa})] \nonumber \\
 & = & \mathbb{E}_{\vec{t},\boldsymbol{\Pi}}[\ln p(\vec{c}|\vec{t},\boldsymbol{\Pi})]+\mathbb{E}_{\vec{t},\boldsymbol{\kappa}}[\ln p(\vec{t}|\boldsymbol{\kappa})]+\mathbb{E}_{\boldsymbol{\boldsymbol{\Pi}}}[\ln p(\boldsymbol{\Pi}|\boldsymbol{\alpha}_{0})]+\mathbb{E}_{\boldsymbol{_{\kappa}}}[\ln p(\boldsymbol{\kappa}|\boldsymbol{\nu}_{0})]  \nonumber\\
 &  & -\mathbb{E}_{\vec{t},\boldsymbol{\Pi},\boldsymbol{\kappa}}[\ln q(\vec{t})]-\mathbb{E}_{\vec{t},\boldsymbol{\Pi}}[\ln q(\boldsymbol{\Pi})]-\mathbb{E}_{\vec{t},\boldsymbol{\kappa}}[\ln q(\boldsymbol{\kappa})] 
\end{eqnarray}

The expectation terms relating to the joint probability of the latent variables, observed variables and the parameters are

\begin{eqnarray}
\mathbb{E}_{\vec{t},\boldsymbol{\Pi}}[\ln p(\vec{c}|\vec{t},\boldsymbol{\Pi})] & = & \sum_{i=1}^{N}\sum_{k=1}^{K}\sum_{j=1}^{J}\mathbb{E}[t_{i}=j]\mathbb{E}\left[\ln\pi_{jc_{i}^{(k)}}^{(k)}\right] \nonumber\\
& = & \sum_{k=1}^{K} \sum_{j=1}^{J} \sum_{l=1}^L N_{jl}^{(k)} \mathbb{E}\left[\ln\pi_{jc_{i}^{(k)}}^{(k)}\right]\\
\mathbb{E}_{\vec{t},\boldsymbol{\kappa}}[\ln p(\vec{t}|\boldsymbol{\kappa})] & = & \sum_{i=1}^{N}\sum_{j=1}^{J}\mathbb{E}[t_{i}=j]\mathbb{E}\left[\ln\kappa_{j}\right] \nonumber\\
& = & \sum_{j=1}^{J}  N_j  \mathbb{E}[\ln\kappa_{j}]\\
\mathbb{E}_{\boldsymbol{\boldsymbol{\Pi}}}[\ln p(\boldsymbol{\Pi}|\mathbf{A}_{0})] & = & \sum_{k=1}^{K}\sum_{j=1}^{J} \left\lbrace-\ln\mathrm{B}\left(\boldsymbol{\alpha}_{0,j}^{(k)}\right)+\sum_{l=1}^{L}\left(\alpha_{0,jl}^{(k)}-1\right)\mathbb{E}\left[\ln\pi_{jl}^{(k)}\right] \right\rbrace \\
\mathbb{E}_{\boldsymbol{_{\kappa}}}[\ln p(\boldsymbol{\kappa}|\boldsymbol{\nu}_{0})] & = & -\ln\mathrm{B}(\boldsymbol{\nu}_{0})+\sum_{j=1}^{J}(\nu_{0,j}-1)\mathbb{E}[\ln\kappa_{j}]
\end{eqnarray}
where $\displaystyle \mathrm{B}(\boldsymbol a) = \frac{ \prod_{l=1}^L \Gamma (a_l)}{\Gamma(\sum_{l=1}^L a_l)}$ is the Beta function and $\Gamma(a)$ is the Gamma function \cite{abramowitz_handbook_1965}. Terms in the lower bound relating to the expectation of the variational distributions $\boldsymbol{q}$ are

\begin{eqnarray}
\mathbb{E}_{\vec{t},\boldsymbol{\Pi},\boldsymbol{\kappa}}[\ln q(\vec{t})] & = & \sum_{i=1}^{N} \sum_{j=1}^{J} \mathbb{E}[t_{i}=j] \ln \mathbb{E}[t_{i}=j] \\
\mathbb{E}_{\vec{t},\boldsymbol{\Pi}}[\ln q(\boldsymbol{\Pi})] & = & \sum_{k=1}^{K}\sum_{j=1}^{J}  \left\lbrace  -\ln \mathrm{B}\left(\boldsymbol{\alpha_{j}^{(k)}}\right)+\sum_{l=1}^{L}\left(\alpha_{j,l}^{(k)}-1\right)\mathbb{E}\left[\ln\pi_{j,l}^{(k)}\right] \right\rbrace \\
\mathbb{E}_{\vec{t},\boldsymbol{\kappa}}[\ln q(\boldsymbol{\kappa})] & = & -\ln\mathrm{B}\left(\mathbf{N}+\boldsymbol{\nu}_{0}\right)+\sum_{j=1}^{J}\left(N_{j}+\nu_{0,j}-1\right)\mathbb{E}[\ln\kappa_{j}]
\end{eqnarray}
where $\mathbf{N} = [N_1, ..., N_J]$ is a vector of counts for each true class. 

Using these equations, the lower bound can be calculated after each pair of \emph{E-step} and \emph{M-step} steps. Once the value of the lower bound stops increasing the algorithm has converged to the optimal approximate solution.

\section{Galaxy Zoo Supernovae} \label{sec:galaxyzooresults}

We tested the model using a dataset obtained from the Galaxy Zoo Supernovae citizen science project \cite{smith_galaxy_2010}. The aim of the project is to
classify candidate supernova images as either ``supernova'' or ``not supernova''. The dataset contains scores given by individual volunteer citizen scientists (base classifiers) to candidates after answering a series of questions. A set of three linked questions are answered by the users, which are hard-coded in the project repository to scores of -1, 1 or 3, corresponding respectively to decisions that the data point is very unlikely to be a supernova, possibly a supernova and very likely a supernova. These scores are our base classifier outputs $\vec c$.

In order to verify the efficacy of our approach and competing methods, we use ``true'' target classifications obtained from full spectroscopic analysis, undertaken as part of the Palomar Transient Factory collaboration \cite{law2009palomar}. We note that this information, is not available to the base classifiers (the users), being obtained retrospectively. We compare IBCC using both variational Bayes (VB-IBCC) and Gibbs sampling (Gibbs-IBCC), using as output the expected values of $t_{i}$. We also tested simple majority voting, weighted majority voting \& weighted sum \cite{littlestone_weighted_2002} and mean user scores, which the Galaxy Zoo Supernovae currently uses to filter results. For majority voting methods we treat both 1 and 3 as a vote for the supernova class.

The complete dataset contains many volunteers that have provided very few classifications, particularly for positive examples, as there are 322 classifications of positive data points compared to 43941 ``not supernova'' examples. Candidate images vary greatly in how difficult they are to classify, so volunteers who have classified small numbers of positive examples may have seen only easy or difficult examples, leading us to infer biased confusion matrices. Including large numbers of volunteers with little data will also affect our inference over true labels and confusion matrices of other decision makers. Therefore, we perform inference over a subsample of the data. Inferred parameters can be used to update the hyperparameters before running the algorithm again over other data point. To infer confusion matrices accurately, we require sufficient numbers of examples for both positive and negative classes. We therefore first select all volunteers that have classified at least 50 examples of each class, then select all data points that have been classified by at least 10 such volunteers; we then include other volunteers that have classified the expected examples. This process produced a data set of 963 examples with decisions produced from 1705 users. We tested the imperfect decision combination methods using five-fold cross validation. The dataset is divided randomly into five partitions, then the algorithm is run five times, each with a different partition designated as the test data and the other partitions used as training data. In the test partition the true labels are withheld from our algorithms and are used only to measure performance. 

\begin{figure}
\centering\subfloat[Receiver operating characteristic (ROC) curves. \label{fig:gz_roc}]{
\includegraphics[width=0.6\textwidth,trim=20mm 0mm 15mm 0mm]{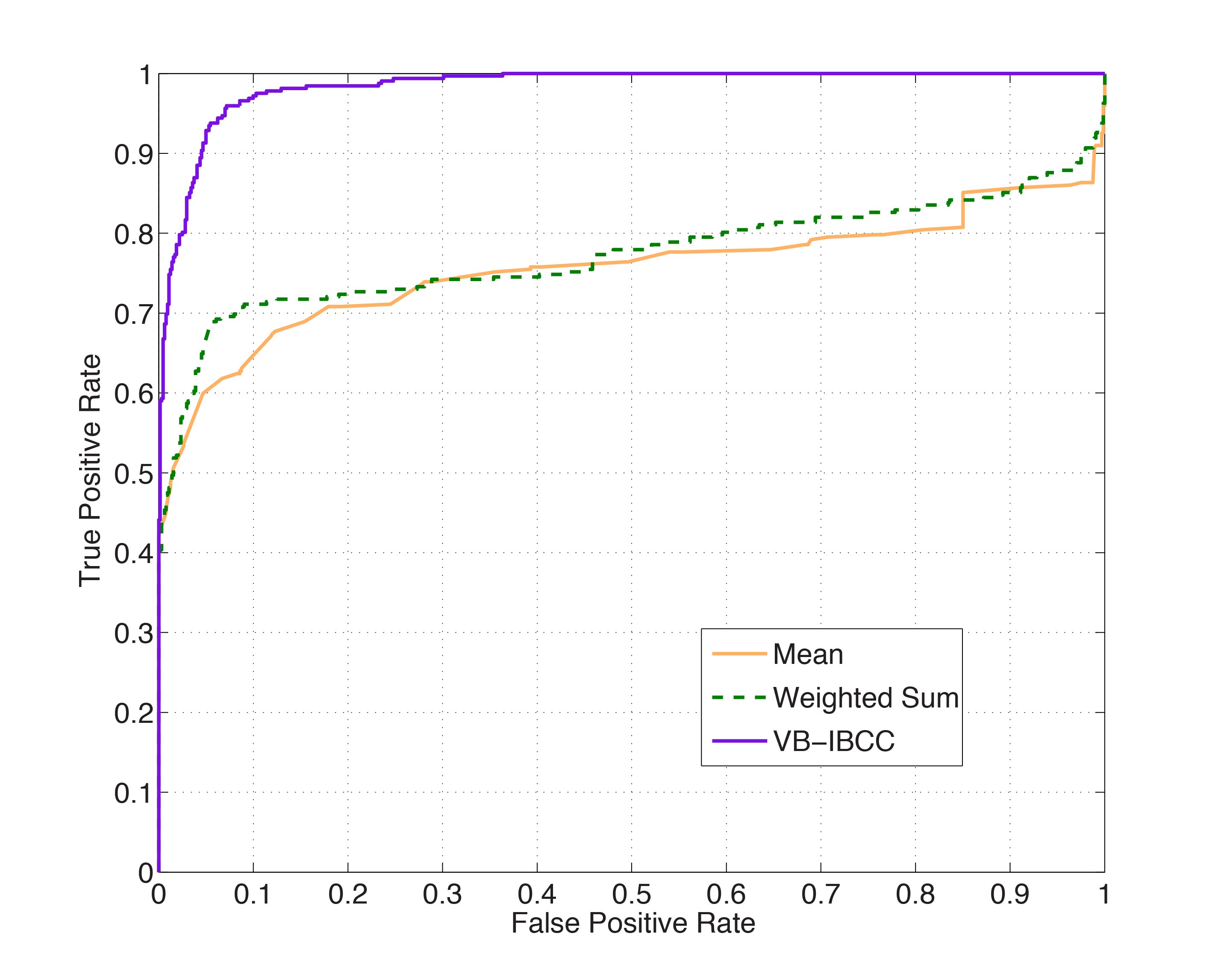}}
\subfloat[Area under the ROC curves (AUCs).\label{fig:gz_auc}] {
\raisebox{15mm}{\begin{tabular}[b]{lr}
\toprule Method & AUC \\
 \midrule \midrule Mean of & \\
Scores & 0.7543 \\
 \midrule Weighted & \\
Sum & 0.7722 \\
 \midrule Simple & \\
Majority  & 0.7809 \\
 \midrule Weighted & \\
 Majority & 0.7378 \\ 
 \midrule Gibbs- & \\
IBCC & 0.9127 \\
 \midrule VB- & \\
IBCC & 0.9840 \\
 \bottomrule
\end{tabular}}}
\caption{Galaxy Zoo Supernovae: ROC curves and AUCs with 5-fold cross validation.}
\end{figure}

Figure \ref{fig:gz_roc} shows the average \emph{Receiver-Operating Characteristic} (ROC) curves \cite{fawcett_introduction_2006} taken across all cross-validation datasets for the mean score, weighted sum and VB-IBCC. Each point on the ROC curve corresponds to a different threshold value; classifier output values above a given threshold are taken as positive classifications and those below as negative. At each threshold value we calculate a true positive rate -- the fraction of positive candidate images correctly identified -- and a false positive rate -- the fraction of negative candidates incorrectly classified as positive. 

The ROC curve for VB-IBCC clearly outperforms the mean of scores by a large margin. Weighted sum achieves a slight improvement on the mean by learning to discount base classifiers each time they make a mistake. The performance of the majority voting methods and IBCC using Gibbs sampling is summarised
by the area under the ROC curve (AUC) in Table \ref{fig:gz_auc}. The AUC gives the probability that a randomly chosen positive instance is ranked higher than a randomly chosen negative instance. Majority voting methods only produce one point on the ROC curve between 0 and 1 as they convert the scores to votes (-1 becomes a negative vote, 1 and 3 become positive) and produce binary outputs. These methods have similar results to the mean score approach, with the weighted version performing slightly worse, perhaps because too much information is lost when converting scores to votes to be able to learn base classifier weights correctly. 

With Gibbs-sampling IBCC we collected samples until the mean of the sample label values converged. Convergence was assumed when the total absolute difference between mean sample labels of successive iterations did not exceed 0.01 for 20 iterations. The mean time taken to run VB-IBCC to convergence was 13 seconds, while for Gibbs sampling IBCC it was 349 seconds. As well as executing significantly faster, VB produces a better AUC than Gibbs sampling with this dataset. Gibbs sampling was run to thousands of iterations with no change in performance observed. Hence it is likely that the better performance of the approximate variational Bayes results from the nature of this dataset; Gibbs sampling may provide better results with other applications but suffers from higher computational costs.

\section{Communities of Decision Makers Based on Confusion Matrices ($\pi$ Communities)} \label{sec:picomms}

In this section we apply a recent community detection methodology to the problem of determining most likely groupings of base classifiers, the imperfect decision makers. Grouping decision makers allows us to observe the behaviours present in our pool of base classifiers and could influence how we allocate classification tasks or train base classifiers. Community detection is the process of clustering  a ``similarity'' or ``interaction'' network, so that classifiers within a given group are more strongly connected to each other than the rest of the graph. Identifying overlapping communities in networks is a challenging task. In recent work \cite{psorakis_overlapping_2011} we presented a novel approach to community detection that infers such latent groups in the network by treating communities
as explanatory latent variables for the observed connections between nodes, so that the stronger the similarity between two decision makers, the more likely it is that they belong to the same community. Such latent grouping is extracted by an appropriate factorisation of the connectivity matrix, where the effective inner rank (number of communities) is inferred by placing shrinkage priors \cite{tan_automatic_2009} on the elements of the factor matrices. The scheme has the advantage of soft-partitioning solutions, assignment of node participation scores to communities, an intuitive foundation and computational efficiency. 

We apply the approach described in \cite{psorakis_overlapping_2011} to a similarity matrix calculated over all the citizen scientists in our study, based upon the expected values of each users' confusion matrix. Expectations are taken over the distributions of the confusion matrices inferred using the variational Bayesian method in Section \ref{sec:VBIBCC} and characterise the behaviour of the base classifiers. Denoting $\mathbb{E}[\pi^{(i)}]$ as the ($3 \times 2$) confusion matrix inferred for user $i$ we may define a simple similarity measure between agents $m$ and $n$ as
\begin{equation}
 V_{m,n} = \exp \left (-\sum_{j=1}^J\mathcal{HD}\left(\mathbb{E}[\boldsymbol{\pi}_j^{(m)}] , \mathbb{E}[\boldsymbol{\pi}_j^{(n)}]\right) \right ), \label{eq:adjacencymatrix}
\end{equation}
where $\mathcal{HD}()$ is the \emph{Hellinger distance} between two distributions, meaning that two agents who have very similar confusion matrices will have high similarity. Since the confusion matrices are multinomial distributions, Hellinger distance is calculated as:
\begin{equation}
\mathcal{HD}\left(\mathbb{E}[\boldsymbol{\pi}^{(m)}_j],\mathbb{E}[\boldsymbol{\pi}^{(n)}_j]\right) = 1 - \sum_{l=1}^L \sqrt{\mathbb{E}[\pi^{(m)}_{jl}] \mathbb{E}[\pi^{(n)}_{jl}]}
\end{equation}
As confusion matrices represent probability distributions, so Hellinger distance is chosen as an established, symmetrical measure of similarity between two probability distributions \cite{bishopbooko}. Taking the exponential of the negative Hellinger distance converts the distance measure to a similarity measure with a maximum of 1, emphasising cases of high similarity.

\begin{figure}
\centering{} {\includegraphics[width=1\textwidth,bb=0 0 480 90]{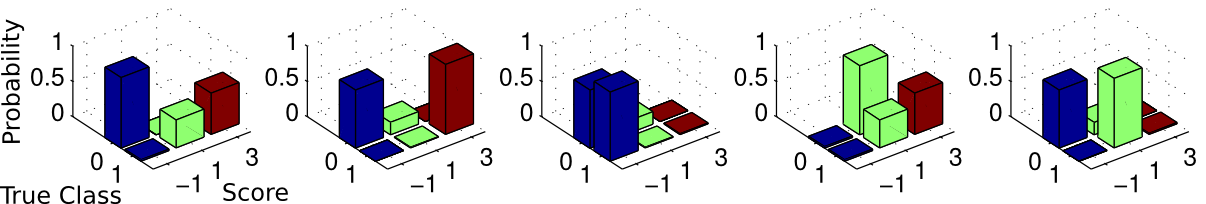}}
\caption{Prototypical confusion matrices for each of the five communities inferred using Bayesian social network analysis. Each graph corresponds to the most central individual in a community, with bar height indicating probability of producing a particular score for a candidate of the given true class.}
\label{fig:confMats}
\end{figure}

Application of Bayesian community detection to the matrix $\boldsymbol{V}$ robustly gave rise to \emph{five} distinct groupings of users. In Figure \ref{fig:confMats} we show the centroid confusion matrices associated with each of these groups of citizen scientists. The centroids are the expected confusion matrices of the individuals with the highest node participation scores for each community. 
The labels indicate the ``true'' class (0 for ``not supernova'' or 1 for ``supernova'') and the preference for the three scores offered to each user by the Zooniverse questions (-1, 1 \& 3). Group 1, for example, indicates users who are clear in their categorisation of ``not supernova'' (a score of -1) but who are less certain regarding the ``possible supernova'' and ``likely supernova'' categories (scores 1 \& 3). Group 2 are ``extremists'' who use little of the middle score, but who confidently (and correctly) use scores of -1 and 3. By contrast group 3 are users who almost always use score -1 (``not supernova'') 
whatever objects they are presented with. Group 4 almost never declare an object as ``not supernova'' (incorrectly) and, finally, group 5 consists of ``non-committal'' users who rarely assign a score of 3 to supernova objects, preferring the middle score (``possible supernova''). It is interesting to note that all five groups have similar numbers of members (several hundred) but clearly each group indicates a very different approach to decision making.

\section{Common Task Communities}\label{sec:taskcomms}

In this section we examine groups of decision makers that have completed classification tasks for similar sets of objects, which we label \emph{common task communities}. Below we outline how these communities and the corresponding confusion matrices could inform the way we allocate tasks and train decision makers. Intelligent task assignment could improve our knowledge of confusion matrices, increase the independence of base classifiers selected for a task, and satisfy human agents who prefer to work on certain types of task. We apply the overlapping community detection method \cite{psorakis_overlapping_2011} to a co-occurrence network for the Galaxy Zoo Supernovae data. Edges connect citizen scientists that have completed a common task, where edge weights $w_{mn}$ reflect the proportion of tasks common to both individuals, such that
\begin{equation}
 w_{mn} = \frac {\mathrm{number\_of\_common\_tasks}(m, n)} {0.5(N^{(m)} + N^{(n)})} \;,
\end{equation}
where $N^{(k)}$ is the total number of observations seen by base classifier $k$. The normalisation term reduces the weight of edges from decision makers that have completed large numbers of tasks, as these would otherwise have very strong links to many others that they proportionally have little similarity to. The edge weights capture the correlation between the tasks that individuals have completed and give the expectation that for a classifier label $c^{(m)}_i$ chosen randomly from our sample, the classifier $n$ will also have produced a label $c^{(n)}_i$. It is possible to place a prior distribution over these weights to provide a fully Bayesian estimate of the probability of classifiers completing the same task. However, this would not affect the results of our community analysis method, which uses single similarity values for each pair of nodes. For decision makers that have made few classifications, edge weights may be poor estimates of similarity and thus introduce noise into the network. We therefore 
filter out decision makers that have made fewer than 10 classifications.

The algorithm found 32 communities for 2131 citizen scientists and produced a strong community structure with modularity of 0.75. Modularity is a measure between -1 and 1 that assumes a strong community structure has more intra-community edges (edges that connect nodes in the same community) than inter-community edges. It is the fraction of intra-community edges minus the expected fraction of intra-community edges for a random graph with the same node degree distribution \cite{girvan_community_2002}. In Galaxy Zoo Supernovae, this very modular community structure may arise through users with similar preferences or times of availability being assigned to the same objects.  Galaxy Zoo Supernovae currently prioritises the oldest objects that currently lack a sufficient number of classifications and assigns these to the next available citizen scientists. It also allows participants to reject tasks if desired. Possible reasons for rejecting a task are that the decision maker finds the task too difficult or 
uninteresting. Common task communities may therefore form where decision makers have similar abilities, preferences for particular tasks (e.g. due to interesting features in an image) or are available to work at similar 
times. When considering the choice of decision makers for a task, these communities could therefore inform who is likely to be available and who will complete a particular 
task.

\begin{figure} [h]
\centering{} {\includegraphics[clip=true,trim=15mm 0mm 0mm 0mm,width=0.4\textwidth]{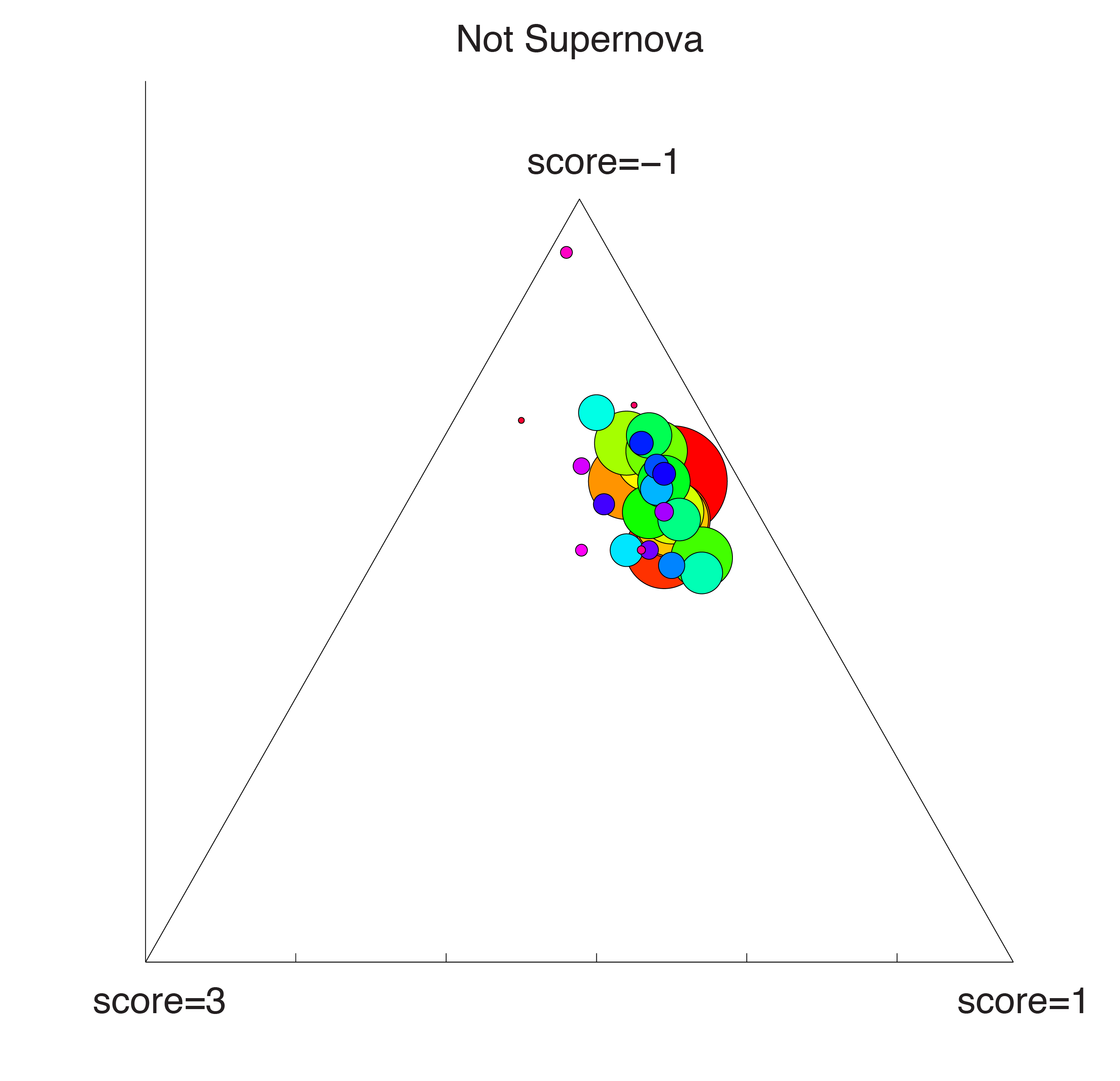}}{\includegraphics[clip=true,trim=15mm 0mm 0mm 0mm,width=0.4\textwidth]{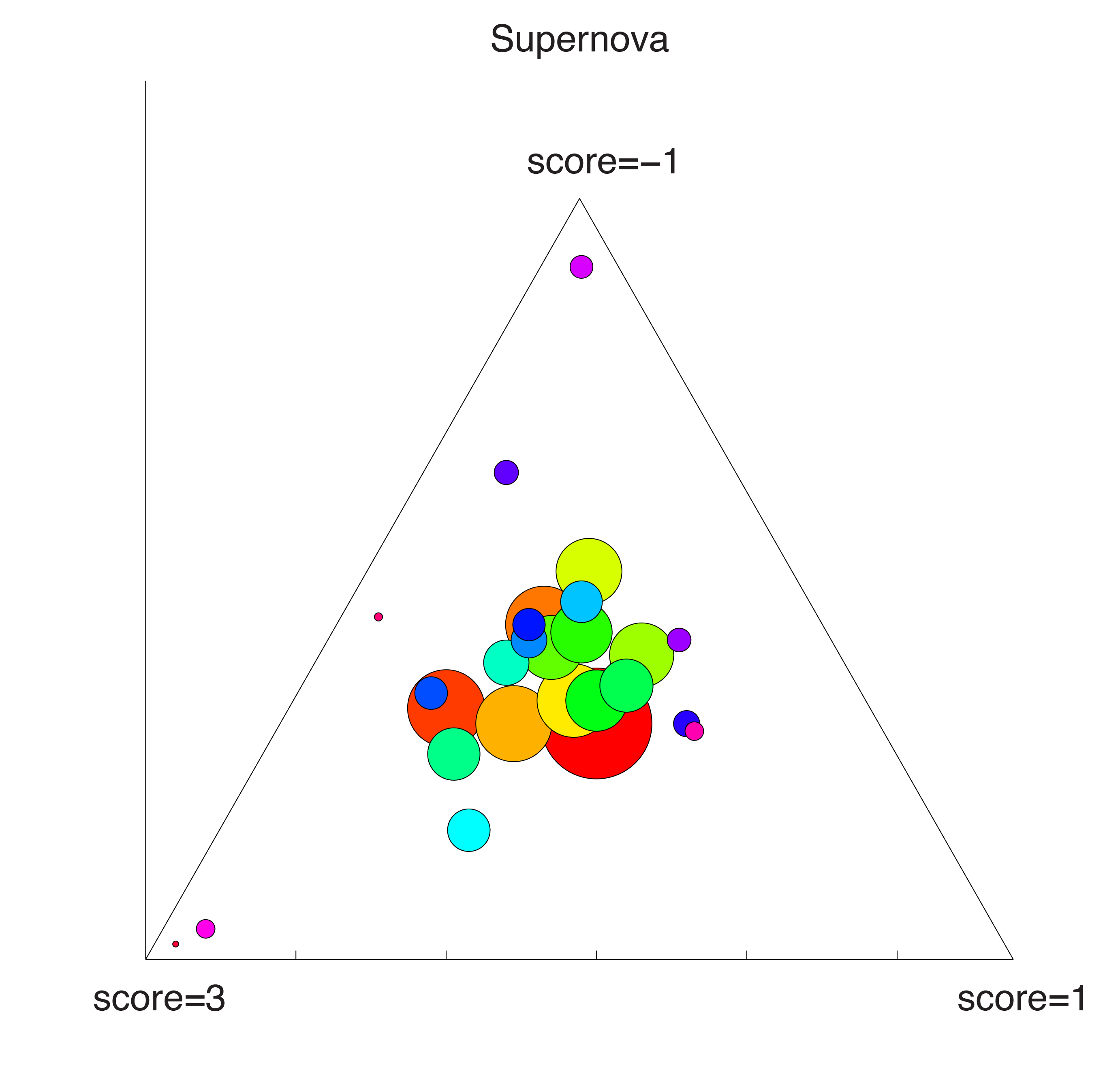}}
\caption{Distribution of means of community members' confusion matrices for all common task communities. Proximity to a vertex indicates the probability of a score given an object with the stated true label class, e.g. in the graph labelled ``Supernova'', a point near the vertex labelled ``score==-1'' indicates a very high probability of giving a decision of -1 when presented with images of a genuine supernova. The left-hand plot shows the mean confusion matrices for $t=0$, i.e. the class ``not a supernova''; the right-hand plot shows the confusion matrices for $t=1$ or ``supernova''. The size of the nodes indicates the number of members of the cluster.}
\label{fig:commTaskTern}
\end{figure}

In Figure \ref{fig:commTaskTern} we plot the distribution of the means of the community members' confusion matrices for each of the true classes. Differences between communities for $t=0$ (not supernova class) are less pronounced than for the $t=1$ (the supernova class). For the latter class we have 5134 observations as opposed to 48791 for the former, so that decision makers see fewer tasks with true label ``supernova''. This means that individual tasks with features that are more likely to elicit a certain base classifier response can have a greater effect on the confusion matrices learned. For instance, some tasks may be easier to classify than others or help a decision maker learn through the experience of completing the task, thus affecting the confusion matrices we infer. As we would expect, some smaller communities have more unusual means as they are more easily influenced by a single community member. The effect of this is demonstrated by the difference between community means in Figure \ref{fig:commTaskTern} for groups of decision makers that have completed common sets of tasks. 

\section{Dynamic Bayesian Classifier Combination} \label{sec:dynibcc}

In real-world applications such as Galaxy Zoo Supernovae, confusion matrices can change over time as the imperfect decision makers learn and modify their behaviour. 
We propose a dynamic variant of IBCC, \emph{DynIBCC}, that models the change each time a decision maker performs a classification task. Using these dynamic confusion matrices we can also observe the effect of each observation on our distribution over a confusion matrix.

\begin{figure}
\centering{} {\includegraphics[width=0.8\textwidth,clip=true,bb=0 70 800 600]{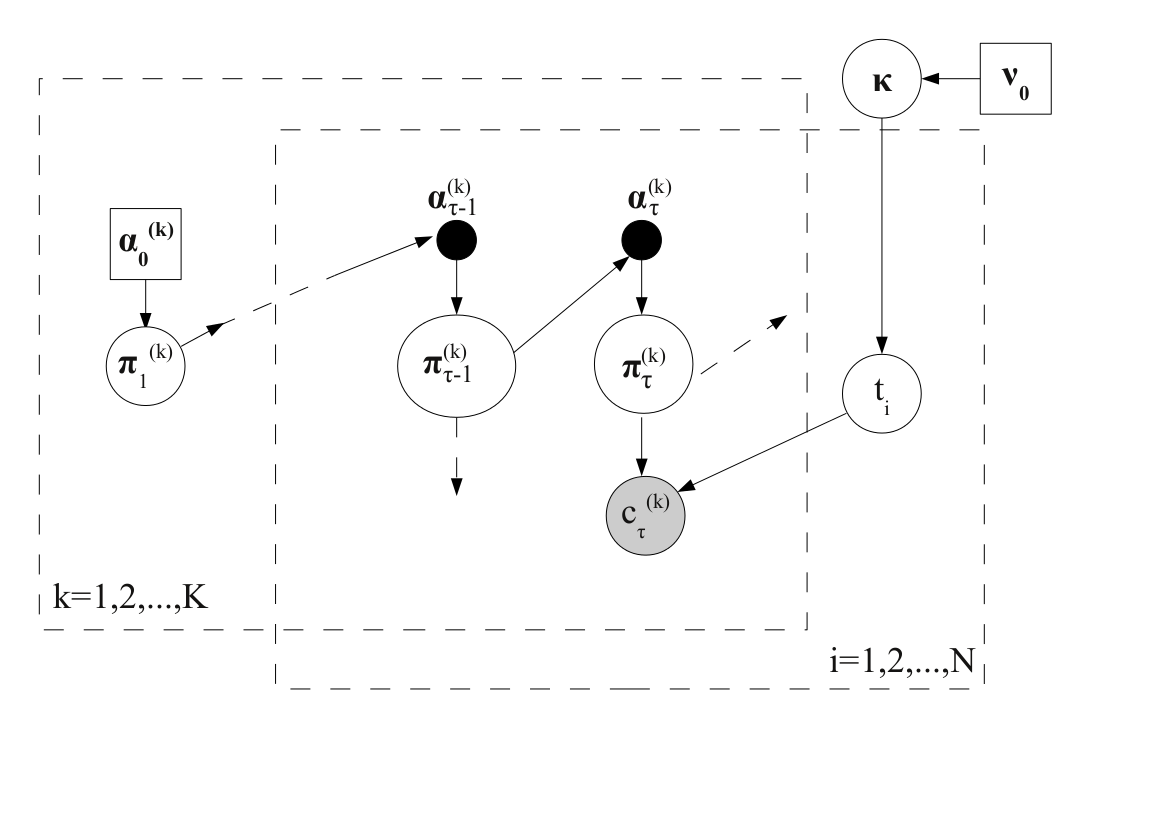}}
\caption{Graphical model for DynIBCC. The dashed arrows indicate dependencies to nodes at previous or subsequent time-steps. Solid black circular nodes are variables calculated deterministically from their predecessors. The shaded node represents observed values, circular nodes are variables with a distribution and square nodes are variables instantiated with point values.}
\label{fig:dynGraphicalModel}
\end{figure}

In DynIBCC, we replace the simple update step for $\boldsymbol{\alpha}_{j}$ given by Equation (\ref{eq:VB_static_alpha_updates}) with an update for every sample classified at time-steps denoted by $\tau$, giving time-dependent parameters $\boldsymbol{\alpha}_{\tau,j}$. Figure \ref{fig:dynGraphicalModel} shows the graphical model for DynIBCC. As we detail in this section, the values of $\boldsymbol{\alpha}_{\tau,j}$ are determined directly (rather than generated from a distribution) from the values of $\boldsymbol{\Pi}$ for the previous and subsequent samples seen by each base classifier $k$. We use a \emph{dynamic generalised linear model} \cite{west_dynamic_1985}, which enables us to iterate through the data updating $\boldsymbol{\alpha}_{\tau}$ depending on the previous value $\boldsymbol{\alpha}_{\tau-1}$. This is the forward pass which operates according to Kalman filter update equations. We then use a Modified Bryson-Frazier smoother \cite{bierman_fixed_1973} to scroll backward through the data, 
updating $\boldsymbol{\alpha}_{\tau}$ based on the subsequent value $\boldsymbol{\alpha}_{\tau+1}$. The backward pass is an extension to the work in \cite{lee_sequential_2010}, where updates are dependent 
only on earlier values of $\boldsymbol{\alpha}$. DynIBCC hence enables us to exploit a fully Bayesian model for dynamic classifier combination, placing distributions over $\boldsymbol{\pi}_{\tau,j}$, while retaining computational tractability by using an approximate method to update the to hyperparameters at each step.

The base classifier $k$ may not classify the samples in the order given by their global indexes $i=1,...,N$, so we map global indexes to time-steps $\tau=1\mathrm,...,T^{(k)}$ using
\begin{equation}
\tau_i^{(k)}=\boldsymbol{f}^{(k)}(i)
 \label{eq:f_of_k}
\end{equation}
The mapping $\boldsymbol{f}^{(k)}$ records the order that $k$ classified items, with time-step $\tau$ being the time-step that sample $i$ was classified.  For an object $i_{unseen}$ not classified by $k$, $\boldsymbol{f}^{(k)}(i_{unseen})=\emptyset$. The inverse of $\boldsymbol{f}^{(k)}(\tau)$ is $i^{(k)}_{\tau}=\boldsymbol{f}^{-1(k)}(\tau)$, a mapping from the time-step $\tau$  $i$ to the object $i$ that was classified at that time-step.

The dynamic generalised linear model allows us to estimate the probability of a classifier producing output $l$ given true label $t_{i_\tau}$:
\begin{equation}
 \tilde{\pi}_{\tau,l} = \pi_{\tau,t_{i_\tau}l} = p\left(c_{i_\tau}=l|t_{i_\tau}\right)
\end{equation}
 in which we omit the superscripts $^{(k)}$ for clarity. We first specify our \emph{generalised linear model} \cite{lee_sequential_2010} by defining a basis function model with the form
\begin{equation}
\tilde{\pi}_{\tau,l}=g\left(\mathbf{h}_{\tau}^{\mathrm{T}}\mathbf{w}_{\tau,l}\right)
\end{equation}
where $\mathbf h_{\tau}$ is a binary input vector of size $J$, with $h_{\tau,t_{i_\tau}} = 1$ and all other values equal to zero, i.e. a binary vector representation of $t_{i_\tau}=j$. The function $g(.)$ is an activation function that maps the linear predictor $\eta_{\tau,l}=\mathbf{h}_{\tau}^{\mathrm{T}}\mathbf{w}_{\tau,l}$ to the probability $\tilde{\pi}_{\tau,l}$ of base classifier response $l$. If we consider each possible classifier output $l$ separately, the value $\tilde{\pi}_{\tau,l}$ is the probability of producing output $l$ and can be seen as the parameter to a binomial distribution. Therefore $g(.)$ is the logistic function, hence
\begin{equation}
\tilde{\pi}_{\tau,l}=\frac{\exp(\eta_{\tau,l})}{1+\exp(\eta_{\tau,l})} \label{eq:logistic}
\end{equation}
 and its inverse, the canonical link function, is the logit function:
\begin{equation}
\eta_{\tau,l}=\mathrm{logit}(\tilde{\pi}_{\tau,l})=\log\left(\frac{\tilde{\pi}_{\tau,l}}{1-\tilde{\pi}_{\tau,l}}\right)
\end{equation}
In the \emph{dynamic} generalised linear model \cite{west_dynamic_1985,lee_sequential_2010}, we track changes to the distribution over $\tilde{\pi}_{\tau,l}$ over time by treating $\mathbf{w}_{\tau,l}$ as a state variable that evolves according to a random walk
\begin{equation}
\mathbf{w}_{\tau,l}=\mathbf{w}_{\tau-1,l}+\mathbf{v}_{\tau,l}\label{eq:randomwalkstate}
\end{equation}
where $\mathbf{v}_{\tau,l}$ is the state noise vector that corresponds to the drift in the state variable over time. We assume that state noise has a distribution where only the first two moments are known, $\mathbf{v}_{\tau,l}\sim(\mathbf{0},q_{\tau,l}\mathbf{I})$, where $\mathbf{I}$ is the identity matrix. The the state noise variance $q_{\tau,l}$ will be estimated from our distributions over $\tilde{\pi}_{\tau,l}$, as explained below.

\subsection{Prior Distributions over the State Variables}

Here we consider Bayesian inference over the state variable $\mathbf{w}_{\tau,l}$. As we place a distribution over $\tilde{\pi}_{\tau,l}$ in our model, we also have a distribution over $\mathbf{w}_{\tau,l}$. Our sequential inference method is approximate since we only estimate the mean and covariance over $\mathbf{w}_{\tau,l}$ rather than the full form of the distribution. At time $\tau$, given observations up to time-step $\tau-1$, the prior state mean at time $\tau$ is $\hat{\mathbf{w}}_{\tau|\tau-1,l}$ and its prior covariance is $\mathbf{P}_{\tau|\tau-1,l}$ . These are related to the posterior mean and covariance from the previous time-step $\tau-1$ by
\begin{eqnarray}
\mathbf{\hat{w}}_{\tau|\tau-1,l} & = & \mathbf{\hat{w}}_{\tau-1|\tau-1,l} \label{eq:what_prior} \\
\mathbf{P}_{\tau|\tau-1,l} & = & \mathbf{P}_{\tau-1|\tau-1,l}+q_{\tau,l}\mathbf{I}. \label{eq:P_prior}
\end{eqnarray} 
We estimate the state noise variance $q_{\tau,l}$ as
\begin{equation}
q_{\tau+1,l} = \max[ u_{\tau|\tau,l}-u_{\tau|\tau-1,l}, 0 ] + z_{\tau,l} \label{eq:q_tau}
\end{equation}
where $u_{\tau|\tau,l}$ and are $u_{\tau|\tau-1,l}$ the variances in the distribution over the classifier outputs $\vec{c}$ after observing data up to time $\tau$ and $\tau-1$ respectively, and $z_{\tau,l}$ is the uncertainty in the classifier outputs. For observations up to time-step $\upsilon$, we define $u_{\tau|\upsilon,l}$ as:
\begin{eqnarray}
u_{\tau|\upsilon,l} & = & \mathbb{V}[\delta_{c_{i_{\tau}}l}|c_{i_1},...,c_{i_{\upsilon}},\vec{t},\boldsymbol{\alpha}_0] \nonumber \\ 
& = & \hat{\pi}_{\tau|\upsilon,l} \left( 1 - \hat{\pi}_{\tau|\upsilon,l} \right)
\end{eqnarray}
where 
\begin{equation}
 \hat{\pi}_{\tau|\upsilon,l} = \mathbb{E}[\tilde{\pi}_{\tau,l}|c_{i_1},...,c_{i_{\upsilon}},\vec{t},\boldsymbol{\alpha}_0] = g\left( \mathbf{h}_{\tau}^{\mathrm{T}} \mathbf{\hat{w}}_{\tau|{\upsilon},l}\right). 
\end{equation}
When the classifier outputs are observed, $z_{\tau,l}$ is zero; when they are not observed, we use $\hat{\pi}_{\tau|\tau-1,l}$ as an estimate of the missing output, so $z_{\tau,l} = u_{\tau|\tau-1,l}$.

From Equations (\ref{eq:what_prior}) and (\ref{eq:P_prior}) we can specify the mean and variance of the prior distribution of $\eta_{\tau,l}$:
\begin{equation}
\hat{\eta}{}_{\tau|\tau-1,l}=\mathbf{h}_{\tau}^{\mathrm{T}}\mathbf{\hat{w}}_{\tau|\tau-1,l}
\end{equation}
\begin{equation}
r_{\tau|\tau-1,l}=\mathbf{h}_{\tau}^{\mathrm{T}}\mathbf{P}_{\tau|\tau-1,l}\mathbf{h}_{\tau}
\end{equation}

We can now use $\hat{\eta}{}_{\tau|\tau-1,l}$ and $r_{\tau|\tau-1,l}$ to estimate the parameters of the prior distribution over $\tilde{\pi}_{\tau,l}$ as follows. The dynamic generalised linear model allows the distribution of the output variable $c_{i_{\tau}}$ to have any exponential family distribution \cite{bishopbooko}. In DynIBCC, the discrete outputs have a multinomial distribution, which is a member of the exponential family, with the Dirichlet distribution as the conjugate prior. Therefore, DynIBCC places a Dirichlet prior over $\boldsymbol{\tilde{\pi}}_{\tau}$ with hyperparameters $\boldsymbol{\tilde{\alpha}}_{\tau}$ that are dependent on the true label $t_{i_\tau}$. If we consider a single classifier output $c_{i_{\tau}}=l$, then
\begin{equation}
 \tilde{\pi}_{\tau,l}  \sim  \mathrm{Beta}(\tilde{\alpha}_{\tau,l}, \beta_{\tau,l}),
\end{equation}
where $\beta_{\tau,l} = \sum_{m=1,l\ne l}^{L}\tilde{\alpha}_{\tau,l}$ and $L$ is the number of possible base classifier output values. 
Since $\boldsymbol{\tilde{\pi}}_{\tau}$ is related to $\eta_{\tau,l}$ by the logistic function (Equation (\ref{eq:logistic})), we can write the full prior distribution over $\eta_{\tau,l}$ in terms of the same hyperparameters:
\begin{eqnarray}
p(\eta_{\tau,l}|c_{i_1},...,c_{i_{\tau-1}},\vec{t}) & = & \frac{1}{\mathrm{B}(\tilde\alpha_{\tau,l},\beta_{\tau,l})} \frac{\exp(\eta_{\tau,l})^{\tilde\alpha_{\tau,l}-1}}{(1+\exp(\eta_{\tau,l}))^{\tilde\alpha_{\tau,l}+\beta_{\tau,l}}}
\end{eqnarray}
where $\mathrm{B}(a,b)$ is the beta function. This distribution is recognised as a beta distribution of the second kind \cite{weatherburn_first_1949}. We can approximate the moments of this prior distribution as follows:
\begin{eqnarray}
\hat{\eta}_{\tau|\tau-1,l}=\mathbb{E}[\eta_{\tau,l}|{c}_{1},...{c}_{t-1},{t}_{i_1},...{t}_{i_\tau}] & = & \Psi(\tilde\alpha_{\tau,l})-\Psi(\beta_{\tau,l})\\
 & \approx & \log\left(\frac{\tilde\alpha_{\tau,l}}{\beta_{\tau,l}}\right)\\
r_{\tau|\tau-1,l}=\mathbb{V}[\eta_{\tau,l}|{c}_{1},...{c}_{t-1},{t}_{i_1},...{t}_{i_\tau}] & = & \Psi'(\tilde\alpha_{\tau,l})+\Psi'(\beta_{\tau,l})\\
 & \approx & \frac{1}{\tilde\alpha_{\tau,l}}+\frac{1}{\beta_{\tau,l}}
\end{eqnarray}
From these approximations we can calculate $\boldsymbol{\tilde\alpha_{\tau}}$ and $\boldsymbol{\beta_{\tau}}$ in terms of $\boldsymbol{\hat{\eta}}_{\tau|\tau-1}$ and $\boldsymbol{r}_{\tau|\tau-1}$:
\begin{eqnarray}
\tilde\alpha_{\tau,l} & = & \frac{1+\exp(\hat{\eta}_{\tau|\tau-1,l})}{r_{\tau|\tau-1,l}}\\ \label{eq:Dyn_priorAlpha}
\beta_{\tau,l} & = & \frac{1+\exp(-\hat{\eta}_{\tau|\tau-1,l})}{r_{\tau|\tau-1,l}}.
\end{eqnarray}
This gives us approximate hyperparameters for the prior distribution over $\boldsymbol{\tilde{\pi}}_{\tau}$.

\subsection{Forward Pass Filtering Steps}

The forward pass filtering steps update the distribution over $\boldsymbol{\tilde{\pi}}_{\tau}$ given an observation of the base classifier output $c_{i_{\tau}}$ at time $\tau$. We calculate the posterior hyperparameters $\tilde\alpha_{\tau|\tau,l}$ and $\beta_{\tau|\tau,l}$ by adding to the prior parameters:
\begin{eqnarray}
\tilde\alpha_{\tau|\tau,l} & = & \tilde\alpha_{\tau,l}+\delta_{c_{i_\tau}l}\\
\beta_{\tau|\tau,l} & = & \beta_{\tau,l}+(1-\delta_{c_{i_\tau}l}).
\end{eqnarray}
In this way we update the pseudo-counts of base classifier output values, as we did in the static IBCC model in Equation (\ref{eq:VB_static_alpha_updates}). The posterior mean and variance are then approximated by 
\begin{eqnarray}
\hat{\eta}_{\tau|\tau,l} & \approx & \log\left(\frac{\tilde\alpha_{\tau|\tau}}{\beta_{\tau|\tau}}\right)\\
r_{\tau|\tau} & \approx & \frac{1}{\alpha_{\tau|\tau}}+\frac{1}{\beta_{\tau|\tau}}.
\end{eqnarray}
Then, we can apply an update to the mean and covariance of the state variable using linear Bayesian estimation, described in \cite{west_bayesian_1997}:
\begin{eqnarray}
\mathbf{\hat{w}}_{\tau|\tau,l} & = & \mathbf{\hat{w}}_{\tau|\tau-1,l}+\mathbf{K}_{\tau,l}\left(\hat{\eta}_{\tau|\tau,l}-\hat{\eta}_{\tau|\tau-1,l}\right)\\
\mathbf{P}_{\tau|\tau,l} & = & \left(\mathbf{I}-\mathbf{K}_{\tau,l}\mathbf{h}_{\tau}^{\mathrm{T}}\right)  \mathbf{P}_{\tau|\tau-1,l}\left(1-\frac{r_{\tau|\tau,l}}{r_{\tau|\tau-1,l}}\right)
\end{eqnarray}
where $\mathbf{K}_{\tau,l}$, the equivalent of the optimal Kalman gain is 
\begin{equation}
\mathbf{K}_{\tau,l}=\frac{\mathbf{P}_{\tau|\tau-1,l}^{\mathrm{T}}\mathbf{h}_{\tau}}{r_{\tau|\tau-1,l}}
\end{equation}
 and $\mathbf{I}$ is the identity matrix. The term $\displaystyle \frac{r_{\tau|\tau,l}}{r_{\tau|\tau-1,l}}$ in the covariance update corresponds to our uncertainty over  $\eta_{\tau,l}$, which we do not observe directly. Linear Bayes estimation gives an optimal estimate when the full distribution over the state variable $\mathbf{w}_{\tau,l}$ is unknown, and therefore differs from a Kalman filter in not specifying a Gaussian distribution over $v_{\tau,l}$ in Equation (\ref{eq:randomwalkstate}). 

To perform the forward pass we iterate through the data: for each time-step we calculate the prior state moments using Equations (\ref{eq:what_prior}) and (\ref{eq:P_prior}), then update these to the posterior state moments $\mathbf{\hat{w}}_{\tau|\tau}$ and $\mathbf{P}_{\tau|\tau}$. The forward pass filtering operates in a sequential manner as posterior state moments from time-step $\tau-1$ are used to calculate the prior moments for the subsequent time-step $\tau$.

\subsection{Backward Pass Smoothing Steps}

After filtering through the data calculating $\mathbf{\hat{w}}_{\tau|\tau,l}$ and $\mathbf{P}_{\tau|\tau,l}$ we then run a backward pass to find the approximate posterior moments given all subsequent data points, $\mathbf{\hat{w}}_{\tau|N,l}$ and $\mathbf{P}_{\tau|N,l}$, and from these the posterior hyperparameters given all data, $\boldsymbol{\tilde\alpha}_{\tau|N}$. The backward pass is a Modified Bryson-Frazier smoother \cite{bierman_fixed_1973}, which updates the distribution using the adjoint state vector $\boldsymbol{\hat{\lambda}}_{\tau,l}$ and adjoint covariance matrix $\boldsymbol{\hat{\Lambda}}_{\tau,l}$ as follows:

\begin{eqnarray}
\mathbf{\hat{w}}_{\tau|N,l} & = & \mathbf{\hat{w}}_{\tau|\tau,l} - \mathbf{P}_{\tau|\tau,l}\boldsymbol{\hat{\lambda}}_{\tau,l}\\
\mathbf{P}_{\tau|N,l} & = & \mathbf{P}_{\tau|\tau,l} - \mathbf{P}_{\tau|\tau,l}\boldsymbol{\hat{\Lambda}}_{\tau,l} \mathbf{P}_{\tau|\tau,l}.
\end{eqnarray}
In our dynamical system the state $\mathbf{w}_{\tau}$ evolves according to Equation (\ref{eq:randomwalkstate}), so $\boldsymbol{\hat{\lambda}}_{\tau}$
and $\boldsymbol{\hat{\Lambda}}_{\tau}$ are defined recursively as the posterior updates from the subsequent step $\tau+1$ given data from $\tau+1$ to $N$. 
\begin{eqnarray}
\boldsymbol{\tilde{\lambda}}_{\tau,l} & = & -\frac{\mathbf{h}_{\tau}} {r_{\tau|\tau-1,l}} \left(\hat{\eta}_{\tau|\tau,l}-\hat{\eta}_{\tau|\tau-1,l}\right) + \left(\mathbf{I}-\mathbf{K}_{\tau,l}\mathbf{h}_{\tau}^{\mathrm{T}}\right)\boldsymbol{\hat{\lambda}}_{\tau,l}\\
\boldsymbol{\hat{\lambda}}_{\tau} & = & \boldsymbol{\tilde{\lambda}}_{\tau+1}\\
\boldsymbol{\hat{\lambda}}_N & = & \mathbf 0 \\
\boldsymbol{\tilde{\Lambda}}_{\tau,l} & = & \frac{\mathbf{h}_{\tau}\mathbf{h}_{\tau}^{\mathrm{T}}}{r_{\tau|\tau-1,1}}\left(1-\frac{r_{\tau|\tau,l}}{r_{\tau|\tau-1,l}}\right) + \left(\mathbf{I}-\mathbf{K}_{\tau,l}\mathbf{h}_{\tau}^{\mathrm{T}}\right)  \boldsymbol{\hat{\Lambda}}_{\tau,l} \left(\mathbf{I}-\mathbf{K}_{\tau,l}\mathbf{h}_{\tau}^{\mathrm{T}}\right)^{\mathrm{T}} \\
\boldsymbol{\hat{\Lambda}}_{\tau} & = & \boldsymbol{\tilde{\Lambda}}_{\tau+1}\\
\boldsymbol{\hat{\Lambda}}_N & = & \mathbf 0
\end{eqnarray}

Estimates for final posterior hyperparameters are therefore given by 
\begin{eqnarray}
\hat{\eta}{}_{\tau|N,l} & = & \mathbf{h}_{\tau}^{\mathrm{T}}\mathbf{\hat{w}}_{\tau|N,l} \\
r_{\tau|N,l} & = & \mathbf{h}_{\tau}^{\mathrm{T}}\mathbf{P}_{\tau|N,l}\mathbf{h}_{\tau} \\
\tilde\alpha_{\tau|N,l} & = & \frac{1+\exp(\hat{\eta}_{\tau|N,l})}{r_{\tau|N,l}}\\ \label{eq:Dyn_postAlpha}
\beta_{\tau|N,l} & = & \frac{1+\exp(-\hat{\eta}_{\tau|N,l})}{r_{\tau|N,l}} \!.
\end{eqnarray}
%

\subsection{Variational Update Equations}

We can now replace the variational distribution for $q^*(\boldsymbol{\pi}_j)$ given by Equation (\ref{eq:VB_qstart_pi}). We continue to omit the $^{(k)}$ notation for clarity. The dynamic model instead uses a variational distribution for each time-step, $q^*(\boldsymbol{\pi}^{(k)}_{\tau,j})$ given by 
\begin{eqnarray}
q^*(\boldsymbol{\pi}_{\tau,j}) & = & \frac{1}{B(\boldsymbol{\alpha}_{\tau|N,j})}\prod_{l=1}^{L}(\pi_{\tau,jl})^{\alpha_{\tau|N,jl}-1}\\
& = & \mathrm{Dir}(\boldsymbol{\pi}_{\tau,j}|\alpha_{\tau|N,j1},...,\alpha_{\tau|N,jL})
\end{eqnarray}
where $\mathrm{Dir()}$ is the Dirichlet distribution with parameters $\boldsymbol{\alpha}_{\tau|N,j}$ calculated according to
\begin{eqnarray}
\alpha_{\tau|N,jl} & = & \frac{1+\exp(\hat{w}_{\tau|N,jl})}{P_{\tau|N,jjl}} \label{eq:DynVB_alpha}
\end{eqnarray}
We calculate $\hat{\mathbf{w}}_{\tau|N}$ and $\mathbf{P}_{\tau|N}$ using the above filtering and smoothing passes, taking the expectation over $\vec t$ so we replace $h_{\tau,j}$ with $\tilde{h}_{\tau,j} = \mathbb{E}_{\vec{t}}[t_{i_\tau} = j]$. Equation (\ref{eq:DynVB_alpha}) is used to derive the hyperparameters for each row of the confusion matrix and thus each possible value of $t_{i_\tau}$; thus it is equivalent to Equation (\ref{eq:Dyn_postAlpha}) with $h_{\tau,j}=1$. This update equation replaces Equation (\ref{eq:VB_static_alpha_updates}) in the static model. The expectation given by Equation (\ref{eq:expectation_lnPi}) becomes
\begin{equation}
\mathbb{E}[\ln \pi_{\tau,jl}] = \Psi(\alpha_{\tau|N,jl}) - \Psi\left(\sum_{m=1}^L \alpha_{\tau|N,jm}\right)
\end{equation}
This can then be used in the variational distribution over ${t_i}$ as follows, replacing Equation (\ref{eq:VB_lnqstar_ti}):
\begin{equation}
 \ln q^*(t_i) = \mathbb{E}[\ln \kappa_{t_i}] + \sum_{k=1}^K \mathbb{E}\left[\ln \pi^{(k)}_{\tau^{(k)}_i,t_i c_i^{(k)}}\right] + \mathrm{const}. \label{eq:DynVB_lnqstar_ti}
\end{equation}

\subsection{DynIBCC Joint and Posterior Distributions}

In DynIBCC, we altered the IBCC model to use the time-dependent confusion matrices, giving the joint distribution over all latent variables and parameters in the model as follows. Here we are considering distributions for all base classifiers and therefore must re-introduce the $k$ superscript notation. 
\begin{eqnarray}
p(\boldsymbol{\kappa},\boldsymbol{\Pi},\vec{t},\vec{c}|\boldsymbol{\alpha}_{0},\boldsymbol{\nu}_{0}) & = & \prod_{i=1}^{N}\left\{\kappa_{t_{i}}\prod_{k\in{\mathbf{C}_i}}\pi_{\tau_{i}^{(k)},t_{i}c_{i}^{(k)}}^{(k)}\right. \nonumber \\ \nopagebreak
&  & \left. p\left(\boldsymbol{\pi}^{(k)}_{\tau^{(k)}_{i}}|t_i,\boldsymbol{\pi}^{(k)}_{\tau^{(k)}_i-1},\boldsymbol{\alpha}^{(k)}_{0,j}\right)\right\}p(\boldsymbol{\kappa}|\boldsymbol{\nu}_{0}) \label{eq:DynIBCCJoint1}
\end{eqnarray}
where $\mathbf{C}_i$ is the set of base classifiers that have completed classification task $i$. There is a change of notation from Equation (\ref{eq:IBCC}) for the static IBCC model, which iterates over all $K$ classifiers. Here we iterate over the set $\mathbf{C}_i$ because $\tau_i^{(k)}$ is undefined for objects that have not been classified by $k$. The static model (Equation \ref{eq:IBCC}) does not have this issue as the confusion matrix is the same for all tasks, and thus Equation (\ref{eq:IBCC}) defines a joint probability over all observed and unobserved base classifier outputs. In DynIBCC, if we wish to determine a distribution over an unobserved base classifier output $c_{i_{unseen}}^{(k)}$ we must also determine a suitable confusion matrix by determining which time-step the unseen classifier output occurred at.

 In Equation (\ref{eq:DynIBCCJoint1}) above, the prior over $\boldsymbol{\pi}^{(k)}_{\tau^{(k)}_{i}}$ can be estimated by finding the mean $\boldsymbol{\bar{\eta}}_{\tau|\tau-1,jl}$ and variance $\bar{r}_{\tau|\tau-1,jl}$ of the linear predictor from its value at the previous time-step, given by:
\begin{equation}
 {\eta}_{\tau-1,jl} = \ln\left(\frac{\pi_{\tau-1,jl}}{\sum_{m=1,\ne l}^L \pi_{\tau-1,jl}}\right)
\end{equation}
The bar notation $\bar\eta$ indicates a variable that is calculated deterministically from the previous state given the value of $\boldsymbol{\pi}^{(k)}_{\tau^{(k)}_{i}}$.Considering the random walk Equation (\ref{eq:randomwalkstate}), where the change $\mathbf{v}_{\tau,l}$ from the previous state has mean $\mathbf{0}$ and covariance $q_{\tau,l}\mathbf{I}$, the moments of the linear predictor are 
\begin{eqnarray}
 {\bar{\eta}}_{\tau|\tau-1,jl} & = & {\eta}_{\tau-1,jl} \\
 \bar{r}_{\tau|\tau-1,jl} & = & q_{\tau,jl}
\end{eqnarray}
where $q_{\tau,l}$ is estimated as per Equations (\ref{eq:q_tau}).
From these values we can calculate the parameters $\bar\alpha_{\tau,jl}$ for a Dirichlet prior over $\boldsymbol{\pi}^{(k)}_{\tau^{(k)}_{i}}$:
\begin{equation}
 p\left(\boldsymbol{\pi}^{(k)}_{\tau^{(k)}_{i}}|t_i,\boldsymbol{\pi}^{(k)}_{\tau^{(k)}_i-1},\boldsymbol{\alpha}^{(k)}_{0,j}\right) =  \mathrm{Dir}\left(\boldsymbol{\pi}_{\tau^{(k)}_{i}}^{(k)}|\boldsymbol{\bar\alpha}^{(k)}_{{\tau^{(k)}_{i}},j1},..., \boldsymbol{\bar\alpha}^{(k)}_{{\tau^{(k)}_{i}},jL} \right).
\end{equation}
For $\tau_i^{(k)}=0$ the parameter $\bar\alpha_{1,jl} = \alpha_{0,jl}$. For $\tau_i^{(k)}>0$ it is given by: 
\begin{equation}
\bar\alpha_{\tau,jl}  =  \frac{1+\exp(\bar{\eta}_{\tau,jl})} {\bar{r}_{\tau|\tau-1,jl}}.
\end{equation}

\subsection{Duplicate Classifications}

The original static model did not allow for duplicate classifications of the same object by the same base classifier. We assumed that even if a base classifier alters their decision when they see an object a second time, the two decisions are likely to be highly correlated and so cannot be treated as independent. However, the dynamic model reflects the possibility that the base classifier may change its own underlying model; therefore responses may be uncorrelated if they are separated by a sufficient number of time-steps or if the confusion matrix changes rapidly over a small number of time-steps. A model that handles dependencies between duplicate classifications at time $\tau_{original}$ and time $\tau_{duplicate}$ may adjust $\boldsymbol{\pi}_{\tau_{original}}^{(k)}$ and $\boldsymbol{\pi}_{\tau_{duplicate}}^{(k)}$ to compensate for correlation. However, in applications where duplicates only occur if they are separated by a large number of time-steps it may be reasonable to treat them as independent 
observations. In cases where duplicates are allowed we index decisions by their time-step as $c^{(k)}_\tau$. For model variants that permit duplicates the joint distribution is hence:

\begin{eqnarray}
\lefteqn{ p(\boldsymbol{\kappa},\boldsymbol{\Pi},\vec{t},\vec{c}|\mathbf{A}_{0},\boldsymbol{\nu}_{0}) } \nonumber\\
& = & p(\boldsymbol{\kappa}|\boldsymbol{v}_{0}) \prod_{i=1}^{N} \kappa_{t_{i}} \! \prod_{k\in{\mathbf{C}_i}} \prod_{\tau{\in}\boldsymbol{f}^{(k)} \!\!(i)} \!\!\! \pi_{\tau,t_{i}c_{\tau}^{(k)}}^{(k)} p\left(\pi_{\tau,t_{i}c_{\tau}^{(k)}}^{(k)}|t_i,\boldsymbol{\pi}^{(k)}_{\tau \! - \! 1},\boldsymbol{\alpha}^{(k)}_{0}\right)\label{eq:DynIBCCJoint2}
\end{eqnarray}

 where as before $\boldsymbol{f}^{(k)}(i)$ maps an object $i$ to the time-step at which $i$ was classified by base classifier $k$. For a sample $i_{unseen}$ not classified by $k$, $\boldsymbol{f}^{(k)}(i_{unseen})=\emptyset$.

\begin{figure}
\centering{} {\includegraphics[width=0.7\textwidth,bb=0 100 700 600]{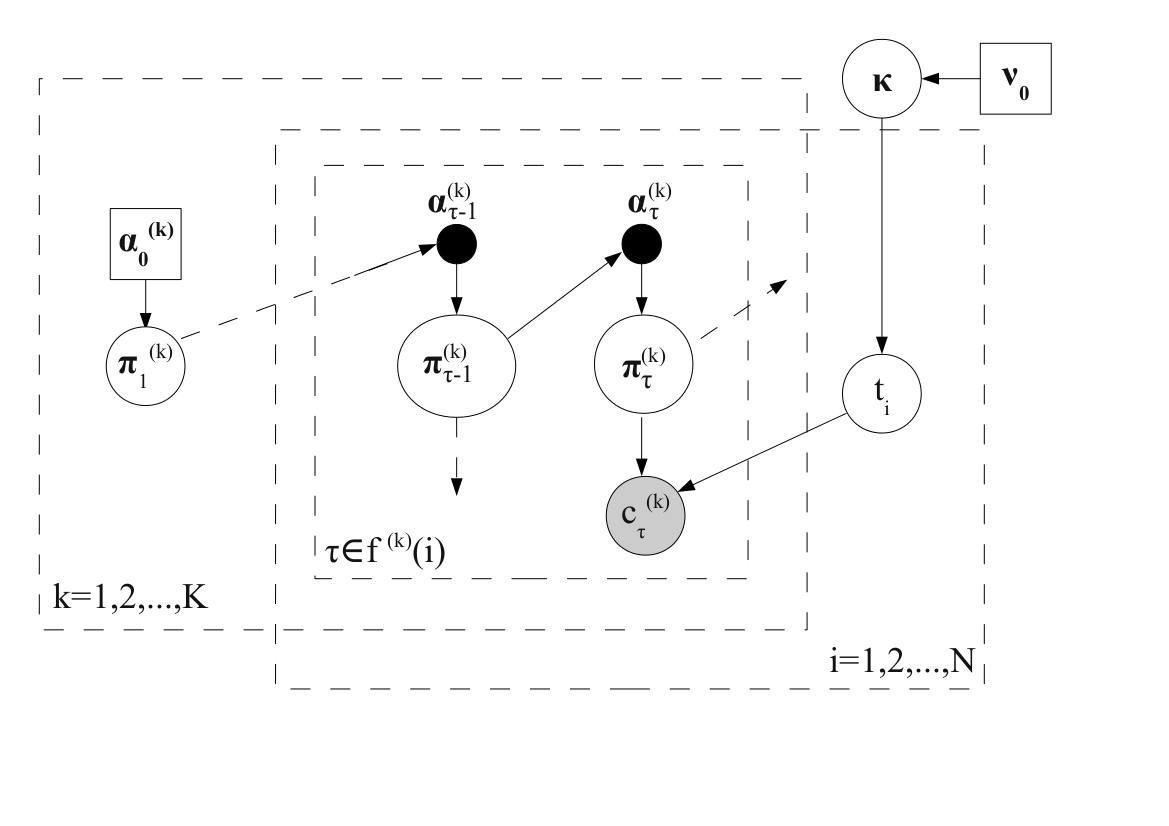}}
\caption{Graphical model for DynIBCC allowing multiple classifications of the same object by the same classifier (duplicate classifications). The dashed arrows indicate dependencies to nodes at previous or subsequent time-steps. Solid black circular nodes are variables calculated deterministically from their predecessors. The shaded node represents observed values, circular nodes are variables with a distribution and square nodes are variables instantiated with point values.}
\label{fig:dynGraphicalModelDupes}
\end{figure}

We must also update Equation (\ref{eq:DynVB_lnqstar_ti}) as to allow duplicates as follows:
\begin{equation}
 \ln q^*(t_i) = \mathbb{E}[\ln \kappa_{t_i}] + \sum_{k\in{\mathbf{C}_i}}\sum_{\tau{\in}\boldsymbol{f}^{(k)}(i)} \mathbb{E}\left[\ln \pi^{(k)}_{\tau,t_i c_i^{(k)}}\right] + \mathrm{const}.
\end{equation}
The resulting graphical model is shown in Figure \ref{fig:dynGraphicalModelDupes}, with an additional plate to allow different time-steps $\tau$ that correspond to the same base classifier $k$ and object $i$.

\subsection{Variational Lower Bound}

We now give the variational lower bound for Dynamic IBCC using the formulation that permits duplicates. We use $\boldsymbol{\Pi}=\left\{ \boldsymbol{\pi}^{(k)}_{\tau,j} | \tau=1,..,T^{(k)},\ j=1,..,J,\ k=1,..,K\right\} $ to refer to the set of confusion matrices for all classifiers and all time-steps.

\begin{eqnarray}
L(q) & = & \int\!\!\! \int\!\!\! \int q(\vec{t},\boldsymbol{\Pi},\boldsymbol{\kappa})\ln\frac{p(\vec{c},\vec{t},\boldsymbol{\Pi},\boldsymbol{\kappa}|\mathbf{A}_{0},\boldsymbol{\nu}_{0})}{q(\vec{t},\boldsymbol{\Pi},\boldsymbol{\kappa})}  \mathrm{d}\vec{t}  \mathrm{d}\boldsymbol{\Pi} \mathrm{d}\boldsymbol{\kappa} \nonumber\\ 
& = & \mathbb{E}_{\vec{t},\boldsymbol{\Pi},\boldsymbol{\kappa}}[\ln p(\vec{c},\vec{t},\boldsymbol{\Pi},\boldsymbol{\kappa}|\mathbf{A}_{0},\boldsymbol{\nu}_{0})]-\mathbb{E}_{\vec{t},\boldsymbol{\Pi},\boldsymbol{\kappa}}[\ln q(\vec{t},\boldsymbol{\Pi},\boldsymbol{\kappa})]  \nonumber\\
& = & \mathbb{E}_{\vec{t},\boldsymbol{\Pi}}[\ln p(\vec{c}|\vec{t},\boldsymbol{\Pi})] + \mathbb{E}_{\boldsymbol{\boldsymbol{\Pi}}}[\ln p(\boldsymbol{\Pi}|\vec{t},\mathbf{A}_{0})] + \mathbb{E}_{\vec{t},\boldsymbol{\kappa}}[\ln p(\vec{t}|\boldsymbol{\kappa})] + \mathbb{E}_{\boldsymbol{_{\kappa}}}[\ln p(\boldsymbol{\kappa}|\boldsymbol{\nu}_{0})]   \nonumber \\
&  & -\mathbb{E}_{\vec{t},\boldsymbol{\Pi},\boldsymbol{\kappa}}[\ln q(\vec{t})]-\mathbb{E}_{\vec{t},\boldsymbol{\Pi}}[\ln q(\boldsymbol{\Pi})]-\mathbb{E}_{\vec{t},\boldsymbol{\kappa}}[\ln q(\boldsymbol{\kappa})]
\end{eqnarray}

The expectation terms relating to the joint probability of the latent variables, observed variables and the parameters are as for the static model in Subsection \ref{sec:staticLowerBound}, except 

\begin{eqnarray}
\mathbb{E}_{\vec{t},\boldsymbol{\Pi}}[\ln p(\vec{c}|\vec{t},\boldsymbol{\Pi})] & = & \sum_{i=1}^{N} \sum_{k\in\mathbf{C}_i} \sum_{\tau{\in}\boldsymbol{f}^{(k)}(i)} \sum_{j=1}^{J} \mathbb{E}[t_{i}=j]\mathbb{E}\left[\ln\pi_{\tau,jc_{\tau}^{(k)}}^{(k)}\right]\\
\mathbb{E}_{\boldsymbol{\Pi}}[\ln p(\boldsymbol{\Pi}|\mathbf{A}_{0})] & = &  \sum_{k=1}^{K}  \sum_{j=1}^{J}  \sum_{\tau=1}^{T^{(k)}} \left\lbrace \!\! -\ln \mathrm{B}\left(\boldsymbol{\bar\alpha}_{\tau,j}^{(k)}\right) \right. \\
& & \left. + \sum_{l=1}^{L}\left(\bar\alpha_{\tau,jl}^{(k)}-1\right)\mathbb{E}\left[\ln\pi_{\tau,jl}^{(k)}\right] \right\rbrace \;\;\;\;\;\;
\end{eqnarray}

In DynIBCC, the expectation over the variational distribution $\boldsymbol{q^*(\pi})$ also differs from static IBCC:

\begin{eqnarray}
\mathbb{E}_{\vec{t},\boldsymbol{\Pi}}[\ln q(\boldsymbol{\Pi})] & = & \sum_{k=1}^{K}  \sum_{j=1}^{J}  \sum_{\tau=1}^{T^{(k)}} \!\! \left\lbrace \!\! - \! \ln\mathrm{B} \! \left(\boldsymbol{\alpha}_{\tau|N,j}^{(k)}\right)   \right. \\
& & \left. + \sum_{l=1}^{L} \left(\alpha_{\tau|N,jl}^{(k)} - 1 \right)\mathbb{E}\left[\ln\pi_{\tau,jl}^{(k)}\right] \right\rbrace 
\end{eqnarray}

\section{Dynamics of Galaxy Zoo Supernovae Contributors} \label{sec:ContributorDynamics}

\begin{figure} [ht!]
\centering{} 
\begin{subfloat}[Volunteer ID 79142\label{fig:piDynamic79142}]{\includegraphics[width=0.4\textwidth,trim=30mm 0mm 10mm 0mm]{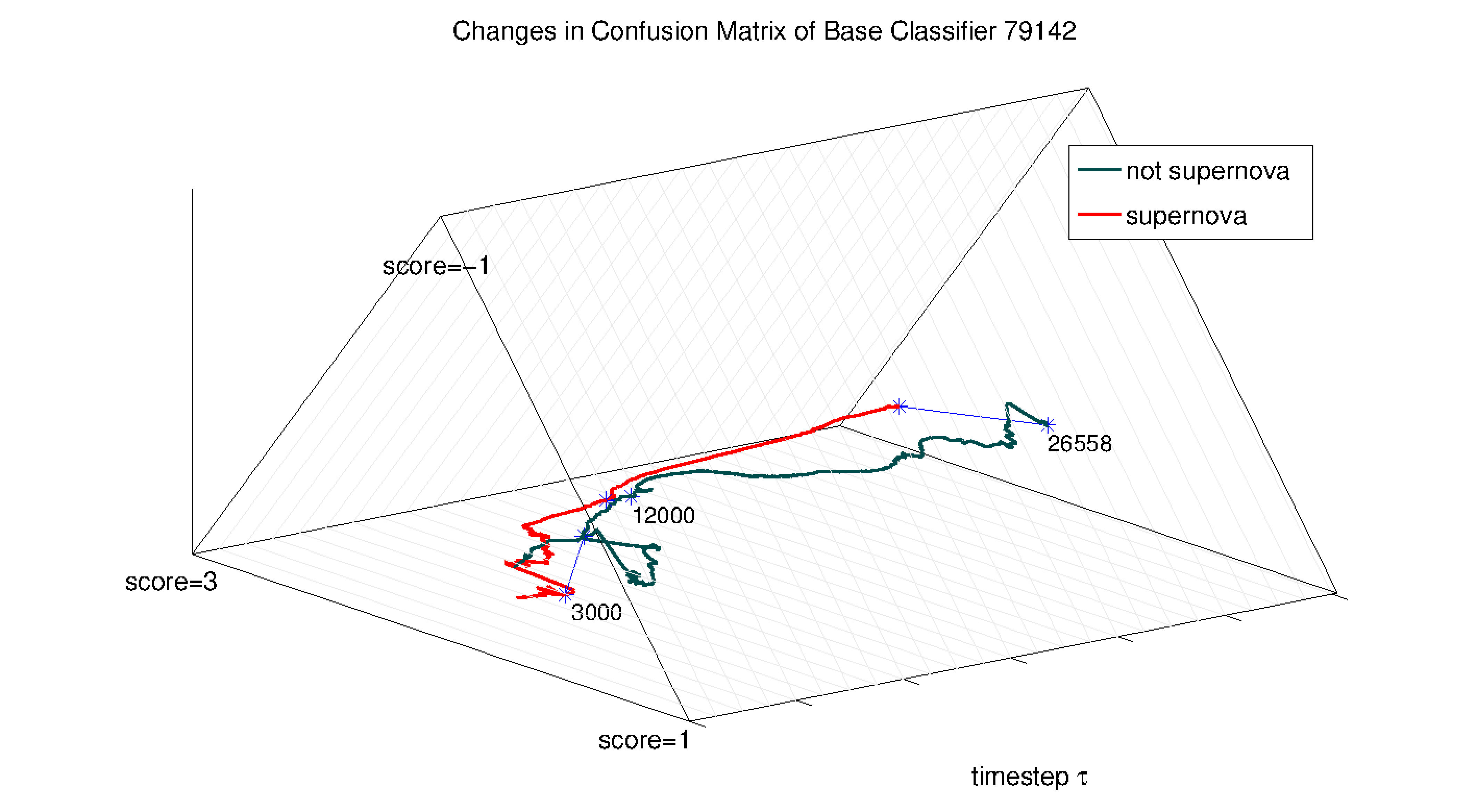}}  \end{subfloat}~~
\subfloat[Volunteer ID 142372\label{fig:piDynamic142372}]{\includegraphics[width=0.4\textwidth,trim=30mm 0mm 10mm 0mm]{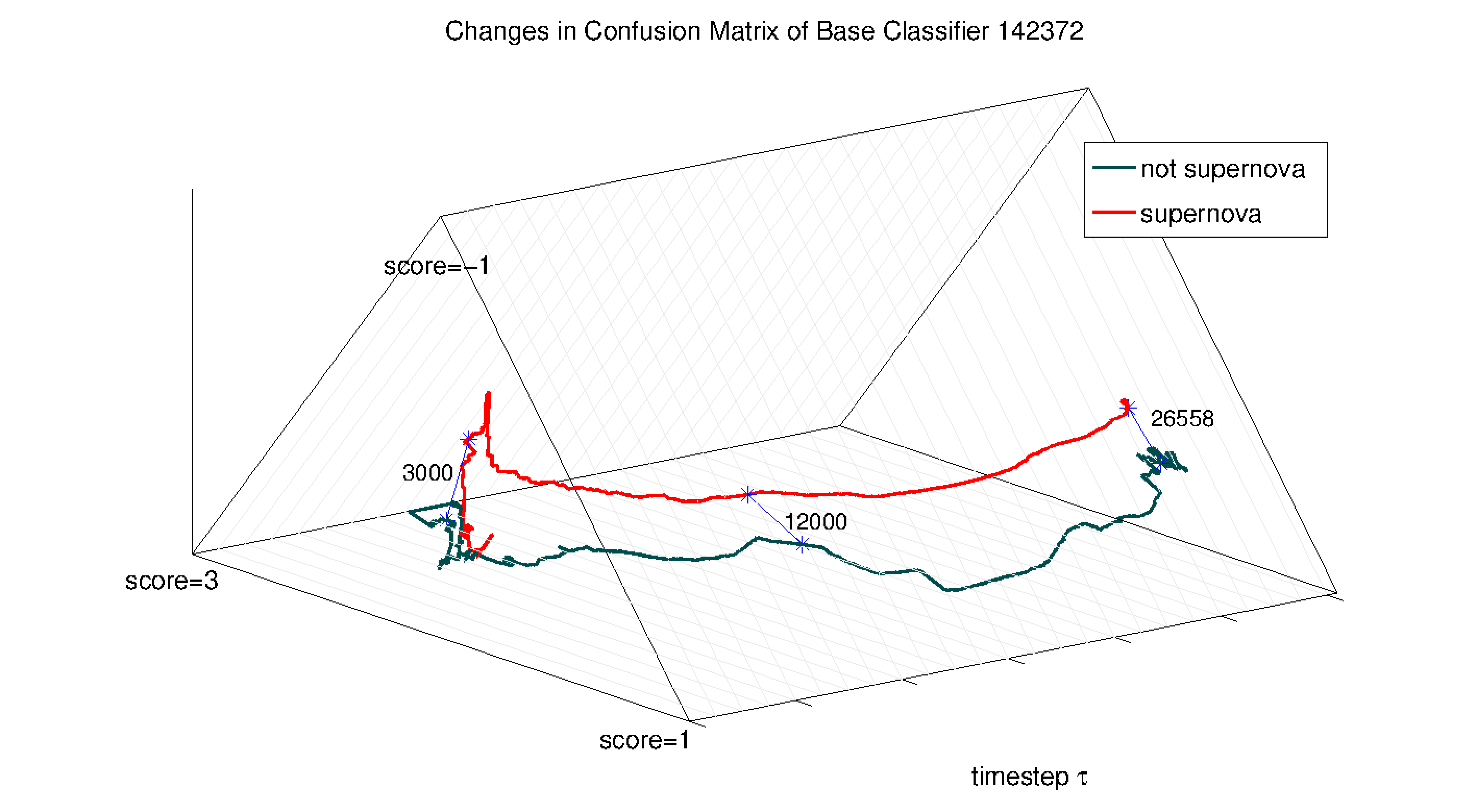}}
\caption{Ternary plot showing the dynamics for Galaxy Zoo Supernovae volunteers. Each line plots the evolution of a row of the confusion matrix corresponding to a particular true class. Proximity of the line to a vertex indicates the probability of generating a certain score for a candidate object with the given true class. Blue '*' markers help the reader align points on the two lines, with a label indicating the global number of observations at that point.}\label{fig:piDynamic1}
\end{figure}

\begin{figure} [ht!]
\centering{} 
\subfloat[Volunteer ID 139963\label{fig:piDynamic139963}]{\includegraphics[width=0.4\textwidth,trim=30mm 0mm 10mm 0mm]{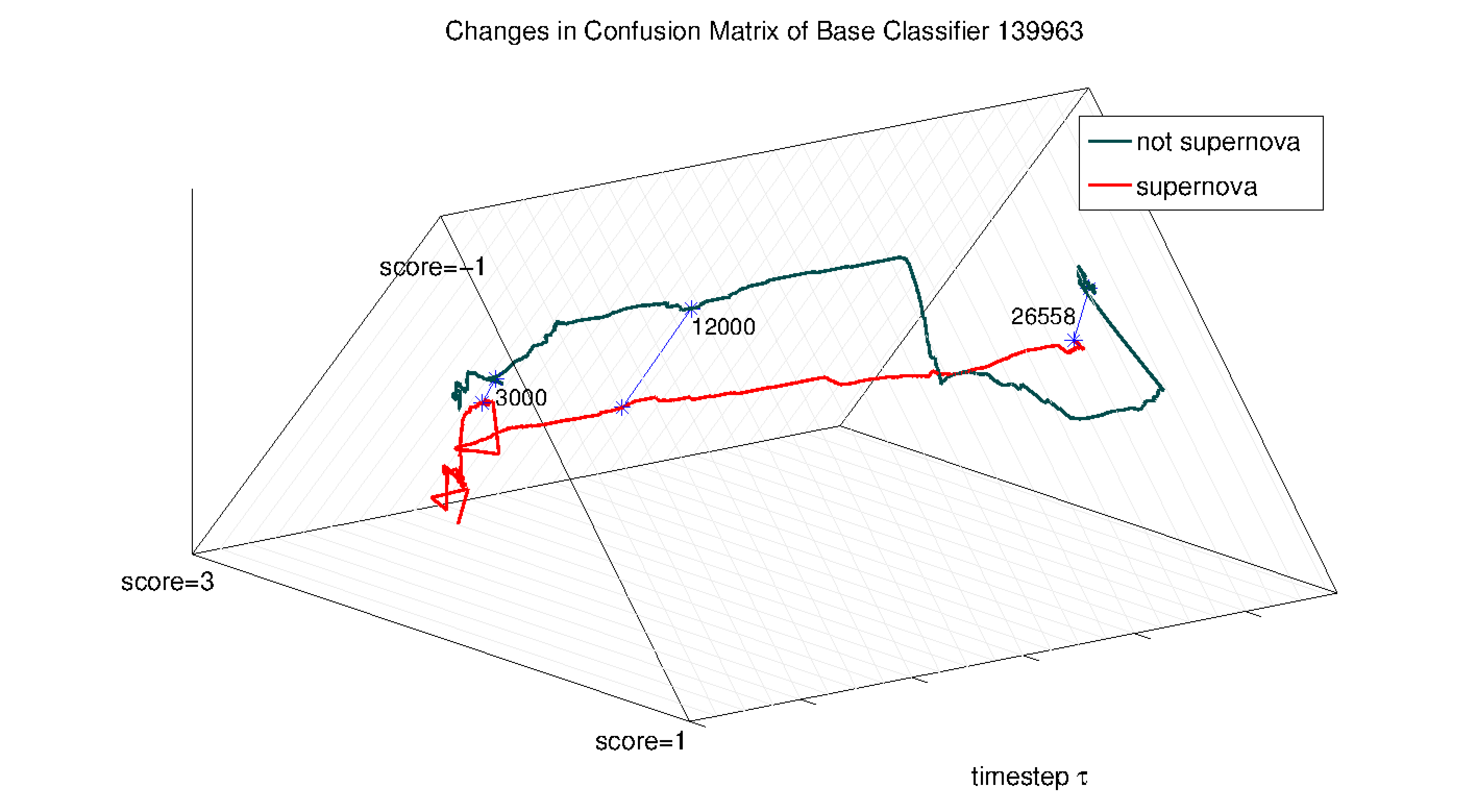}}  ~~
\subfloat[Volunteer ID 259297\label{fig:piDynamic259297}]{\includegraphics[width=0.4\textwidth,trim=30mm 0mm 10mm 0mm]{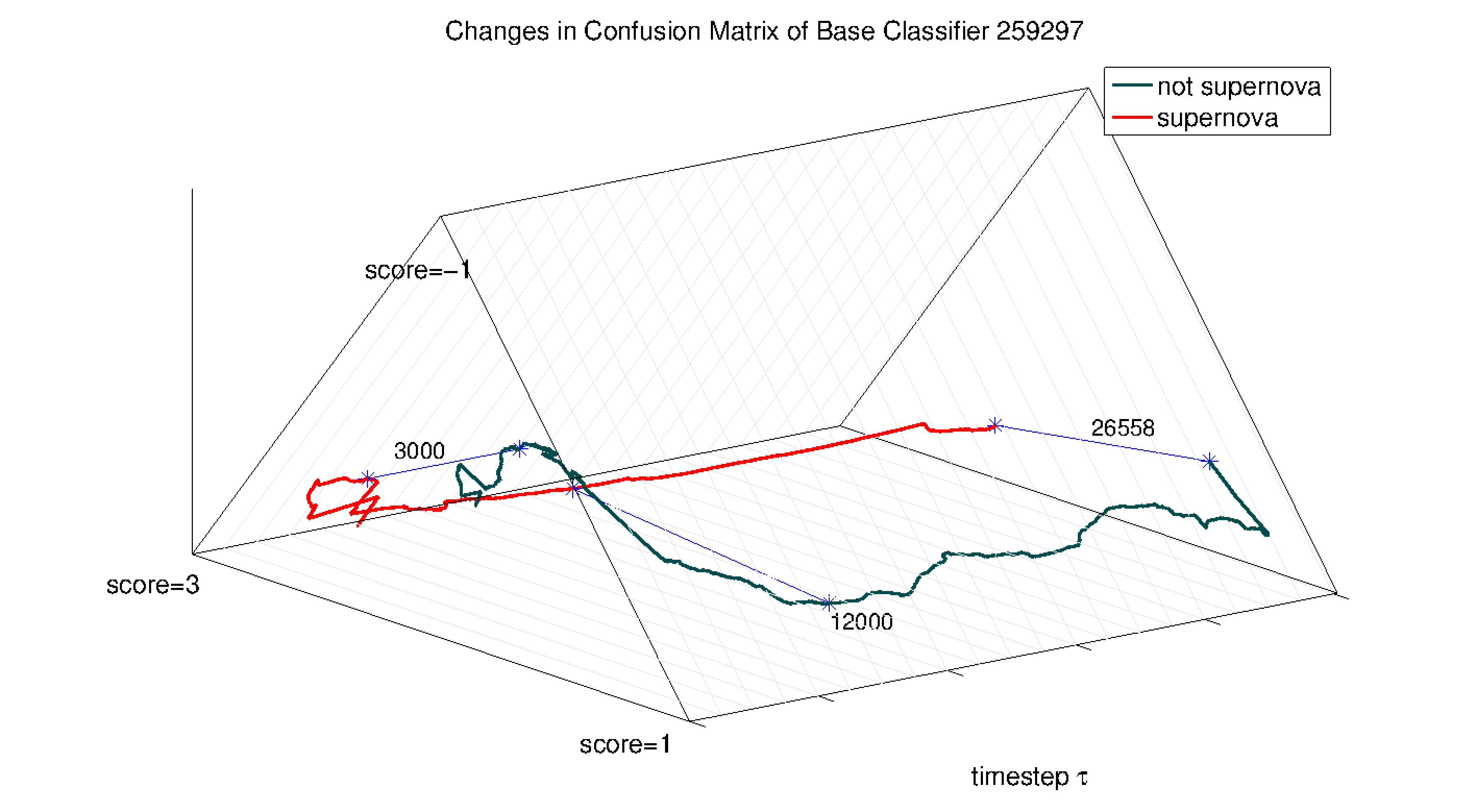}}
\caption{Ternary plot showing the dynamics for Galaxy Zoo Supernovae volunteers. Each line plots the evolution of a row of the confusion matrix corresponding to a particular true class. Proximity of the line to a vertex indicates the probability of generating a certain score for a candidate object with the given true class.Blue '*' markers help the reader align points on the two lines, with a label indicating the global number of observations at that point.}\label{fig:piDynamic2}
\end{figure}

We applied the variational Bayesian DynIBCC to the Galaxy Zoo Supernovae data from Section \ref{sec:galaxyzooresults} to examine the changes to individual confusion matrices. There is a large variation in the dynamics of different decision makers but in many there are sustained drifts in a particular direction. We give examples of the different types of changes found for different base classifiers -- the Galaxy Zoo Supernovae volunteers -- in Figures \ref{fig:piDynamic1} and (\ref{fig:piDynamic2}). These are ternary plots of the expected confusion matrices at each time-step or observation of the decision maker. Each line corresponds to the confusion vector $\boldsymbol \pi_j^{(k)}$ for true class $j$. To help the reader time-align the traces, certain time-steps have been labelled with a blue marker and edge between the two confusion vectors, with the label denoting the global number of observations for all base classifiers at this point. The example volunteers classified 29651, 21933, 23920 and 20869 
candidates respectively. The first example is Figure \ref{fig:piDynamic79142}, which shows a confusion matrix with some drift in the earlier time-steps and some small fluctuations later on for the ``not supernova'' class. The decision maker shown in Figure \ref{fig:piDynamic142372} appears to have a more sustained drift away from scores of 3 for both true classes. The earliest changes in both these decision makers, such as for the ``supernova'' class in Figure \ref{fig:piDynamic142372}, appear to be a move away from the prior, which affects the first data points most. The last two examples show more sudden changes. In Figure \ref{fig:piDynamic139963} we see a very significant change in the later observations for the ``not supernova'' class, after which the confusion matrix returns to a similar point to before. Figure \ref{fig:piDynamic259297} shows little initial change followed by a sudden change for the ``not supernova'' class, which then becomes fairly stable. The dynamics observed were inferred over a 
large number of data points, suggesting that the longer trends are due to genuine changes in performance of base classifiers over time. Smaller fluctuations may be due to bias in the observations (e.g. for a task that is very easy) or a change in behaviour of the citizen scientists, but sustained changes after the initial move away from the priors are more suggestive of a change in behaviour or in the information presented to agents when making classifications. For all four examples, there are more initial fluctuations, which could relate to the way that new volunteers adapt when they complete their first tasks.

\section{Dynamics of $\pi$ Communities} \label{sec:picommsDynamics}

We now apply the community analysis method used in Section \ref{sec:picomms} to the dynamic confusion matrices to examine the development of the community structure over time. After different a number of observations $s$ we run the community detection method \cite{psorakis_overlapping_2011} over an adjacency matrix, using equation  (\ref{eq:adjacencymatrix}) with the most recent confusion matrices for all base classifiers observed up to $s$.

\begin{figure} [!ht]
 \centering{} 
\subfloat[3000 observations.]{\includegraphics[clip=true,trim=0mm 30mm 10mm 40mm,width=0.65\textwidth]{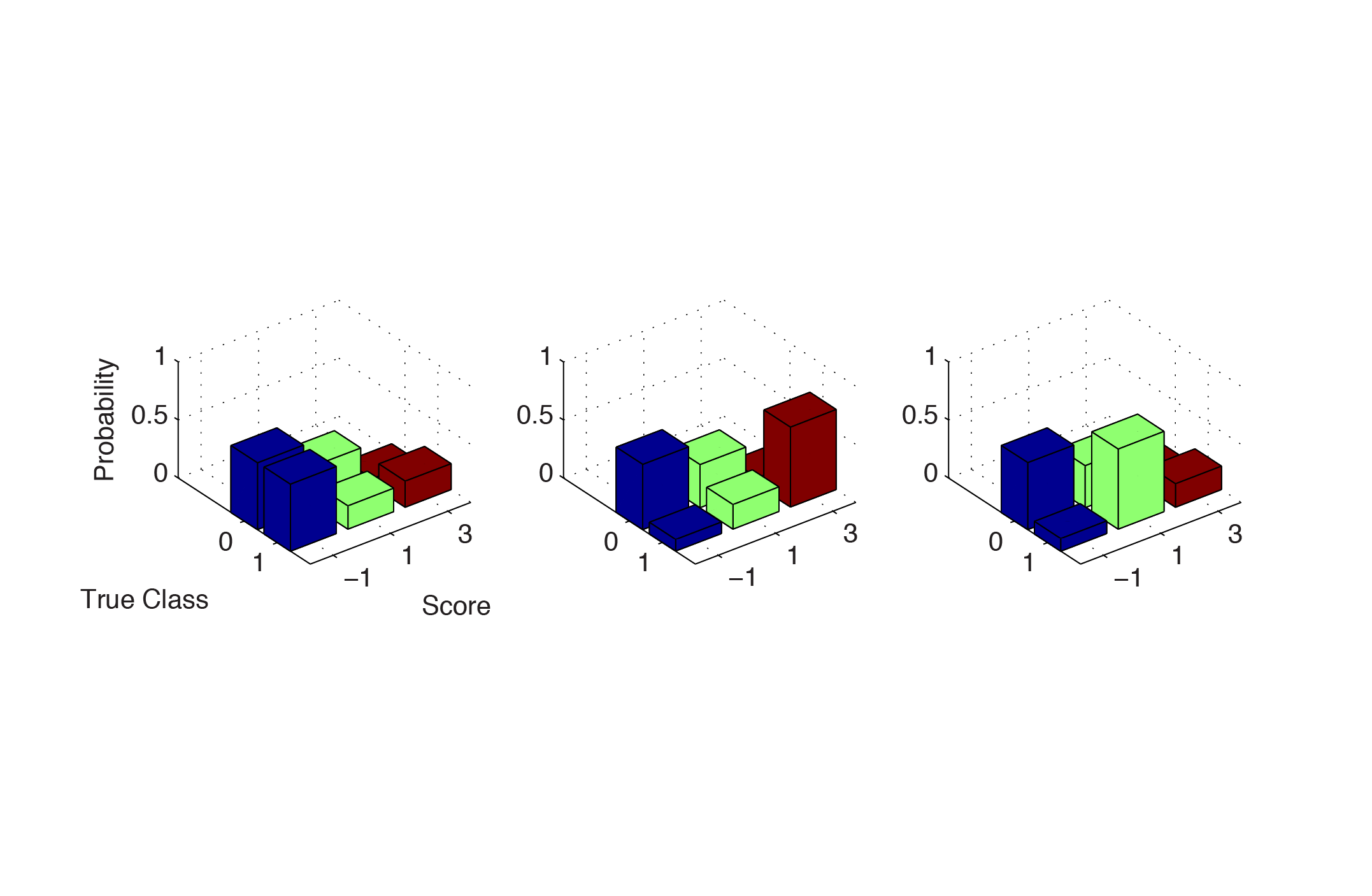}} \\
\subfloat[12000 observations.]{\includegraphics[clip=true,trim=0mm 30mm 10mm 40mm,width=0.9\textwidth]{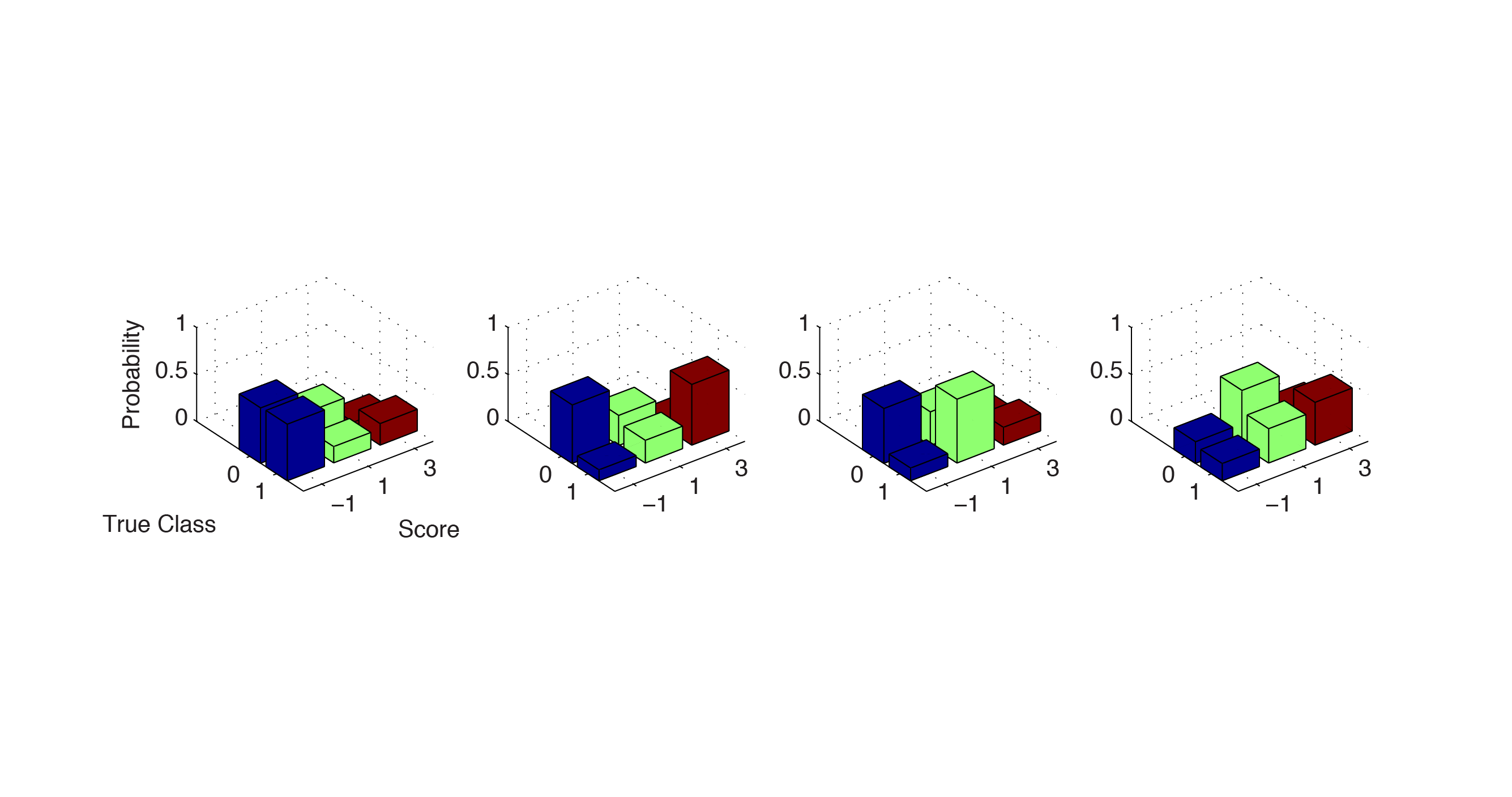}} \\
\subfloat[26558 observations.]{\includegraphics[clip=true,trim=30mm 0mm 30mm 10mm,width=1\textwidth]{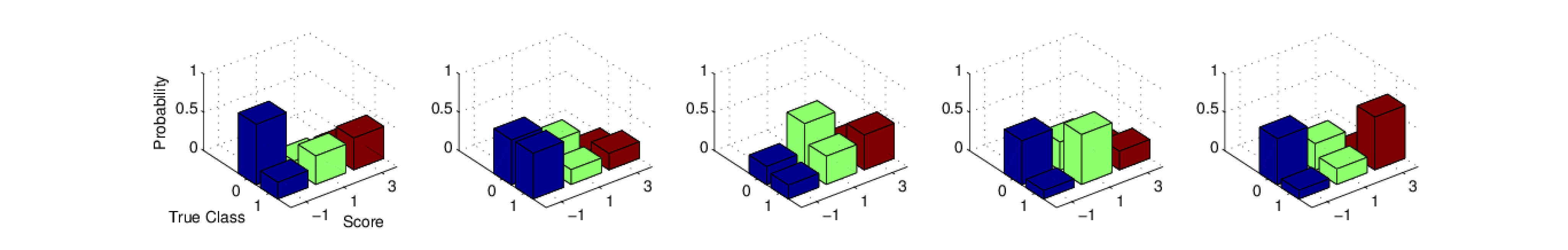}}
\caption{$\boldsymbol\pi$ communities: means over the expected confusion matrices of community members after different numbers of observations. At each time point we see a new community appear while previous communities persist with similar means.}
\label{fig:piCommunityDynamicsMeans}
\end{figure}

In Figure \ref{fig:piCommunityDynamicsMeans} we see how the same communities emerge over time as we saw in Section \ref{sec:picomms} in Figure \ref{fig:confMats}. Initially, only three communities are present, with those corresponding to groups 4 (``optimists'') and 1 (``reasonable'') in Figure \ref{fig:confMats} only appearing after 1200 and 26558 observations. The ``reasonable'' group is the last to emerge and most closely reflects the way the designers of the system intend good decision makers to behave. It may therefore appear as a result of participants learning, or of modifications to the user interface or instructions as the Galaxy Zoo Supernovae application was being developed to encourage this behaviour.

\begin{figure} [!ht]
\centering{}{\includegraphics[clip=true,trim=80mm 100mm 90mm 30mm,width=0.7\textwidth]{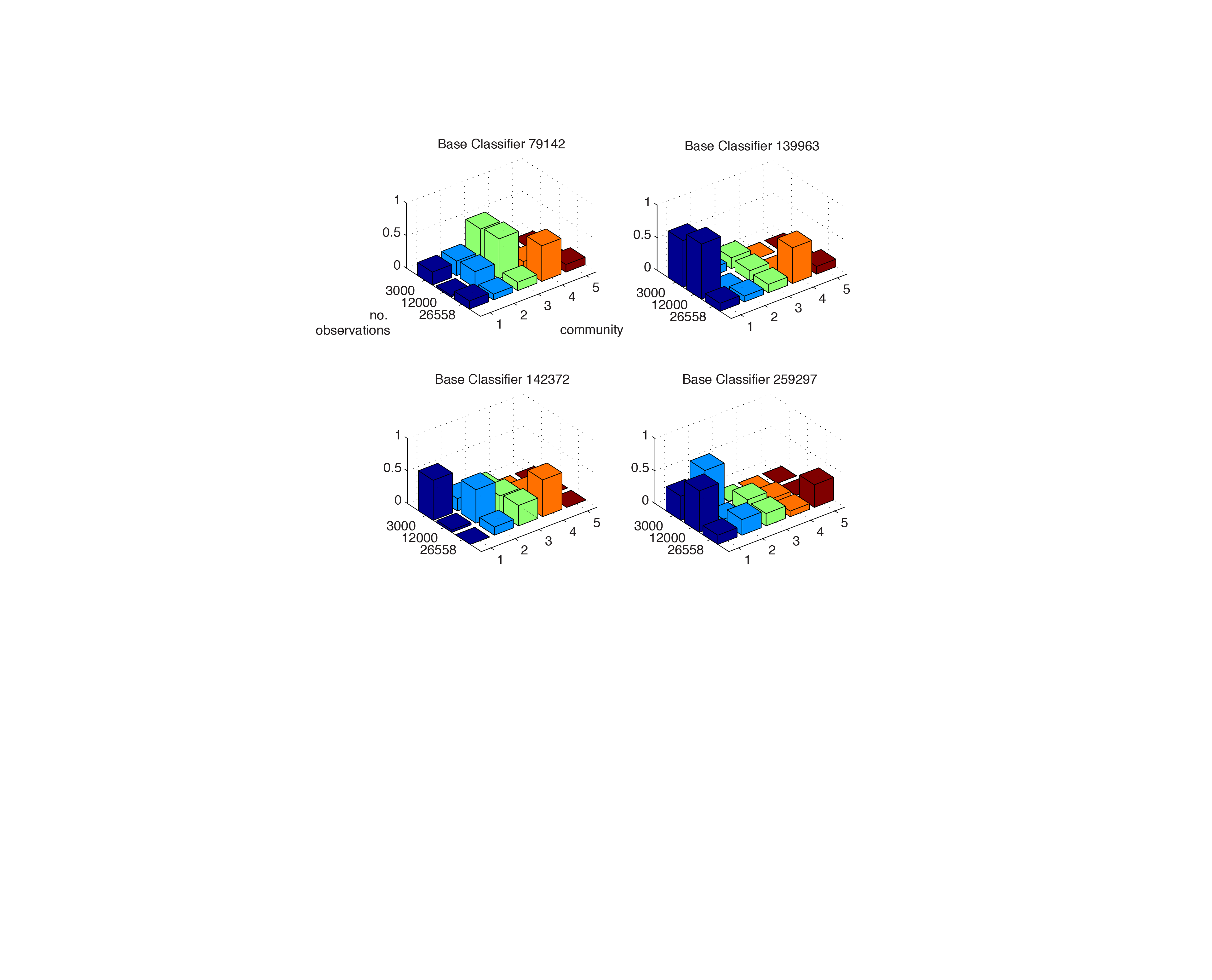}}
\caption{Node participation scores for the $\boldsymbol{\pi}$ communities for selected individuals shown in Section \ref{sec:ContributorDynamics} after different numbers of observations. Each bar corresponds to the individual's participation score in a particular community after running the community analysis over all observations up to that point. Participation scores close to one indicate very strong membership of a community. The node participation scores for one number of observations sum to one over the communities 1 to 5.}
\label{fig:changingmembership}
\end{figure}

We also note that agents switch between communities. In Figure \ref{fig:changingmembership} we show the node participation scores at each number of observations $s$ for the individuals we examined in Section \ref{sec:ContributorDynamics}. Community membership changes after significant changes to the individual's confusion matrix. However, the communities persist despite the movement of members between them.

\section{Dynamics of Common Task Communities} \label{sec:taskcommsDynamics}

\begin{figure} [!ht]
\centering{} 
\includegraphics[clip=true,trim=0mm 0mm 0mm 0mm,width=0.4\textwidth,trim=15mm 0mm 0mm 0mm]{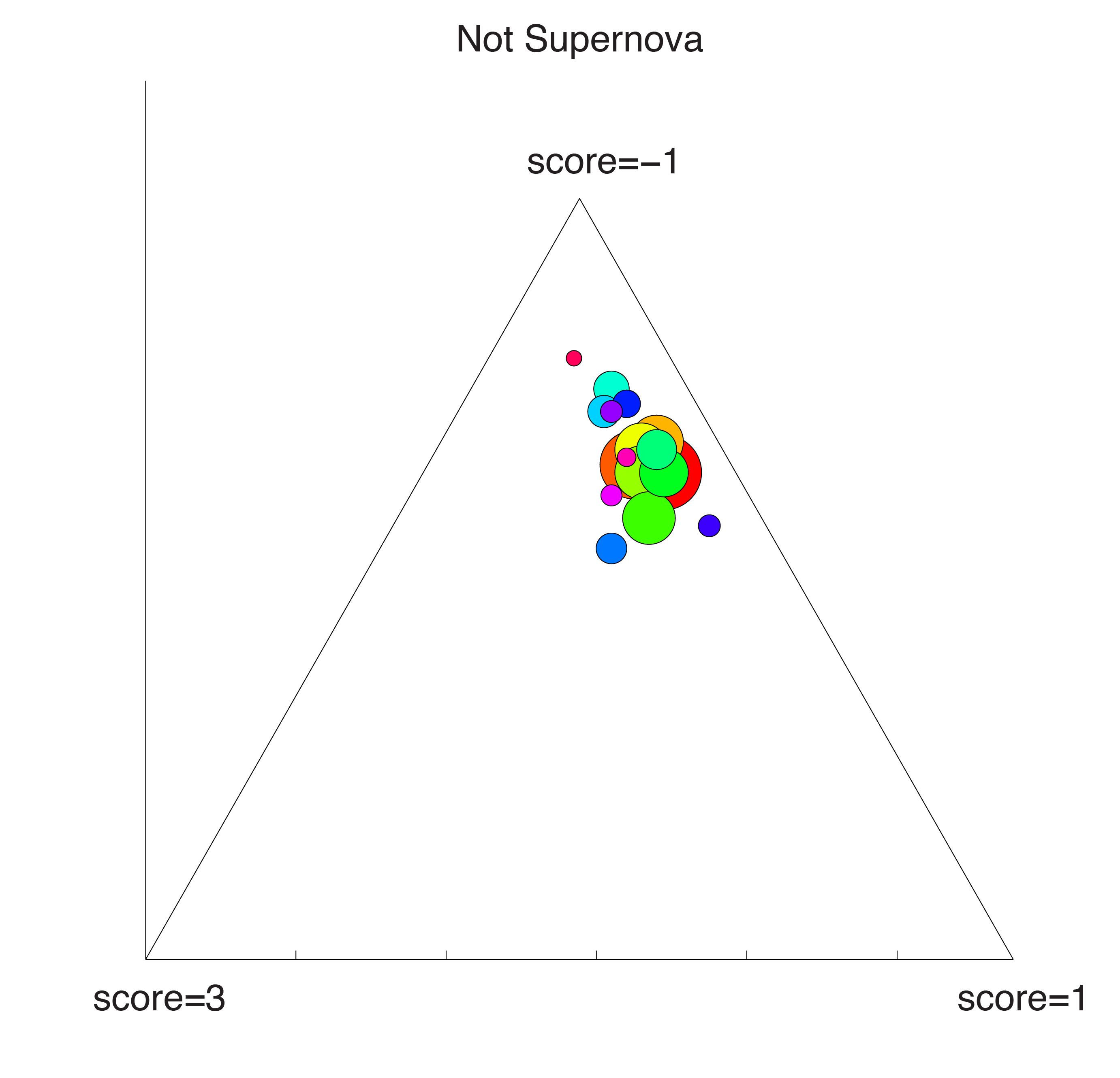}\includegraphics[clip=true,trim=0mm 0mm 0mm 0mm,width=0.4\textwidth,trim=15mm 0mm 0mm 0mm]{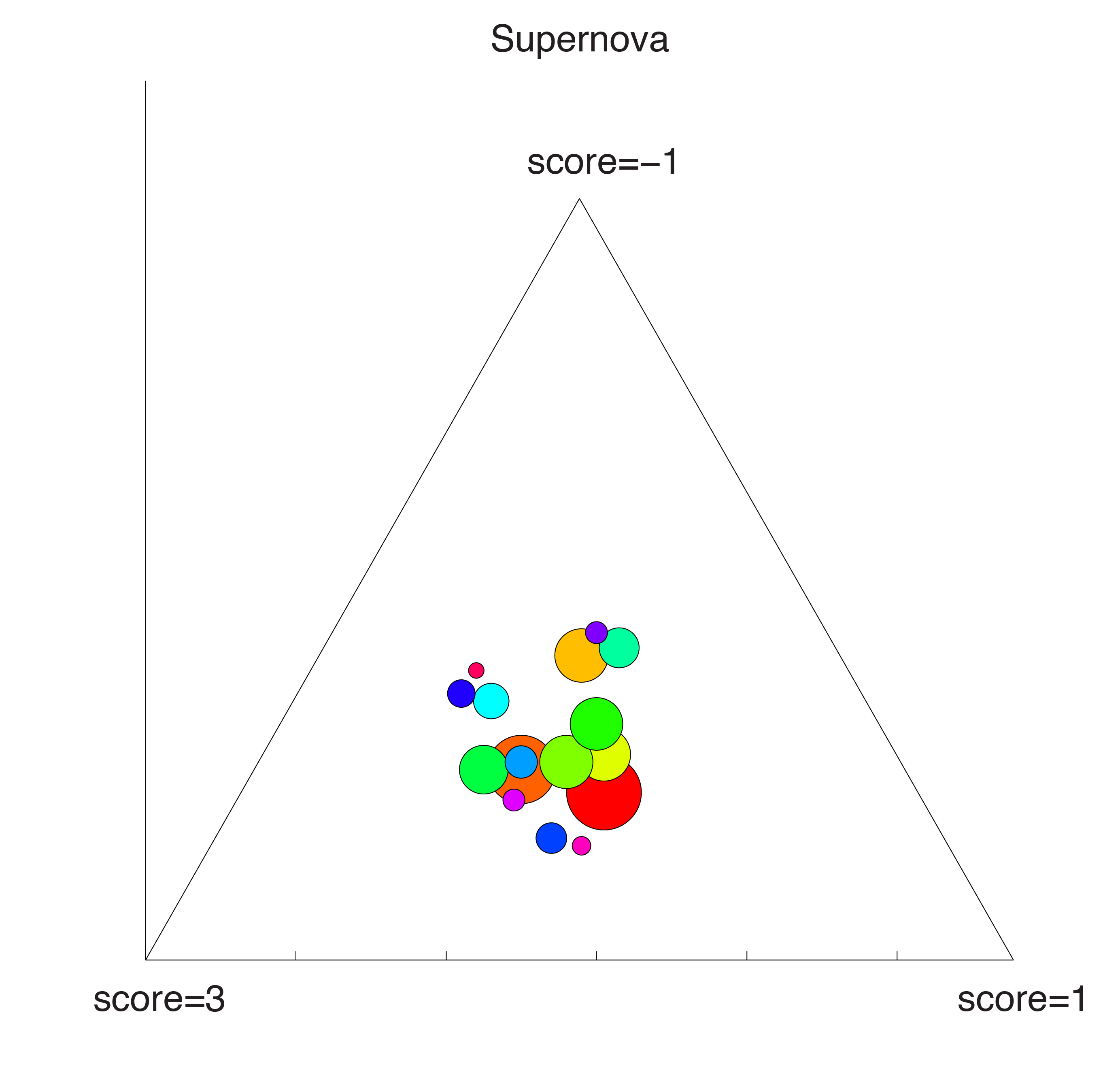}
\caption{Common task communities: ternary plot of means of expected confusion matrices for community members after 50,000 observations. Each point corresponds to one common task community and represents the mean of $\mathbb{E}[\boldsymbol{\pi}_j]$ for the community members. Proximity of a point to a vertex indicates the probability of outputting a particular score when presented with an object of true class ``supernova'' or ``not supernova''.}\label{fig:tasksCommunityDynamicsMeans50000}
\end{figure}

\begin{figure} [!ht]
\centering{} 
 \includegraphics[clip=true,trim=0mm 0mm 0mm 0mm,width=0.4\textwidth,trim=15mm 0mm 0mm 0mm]{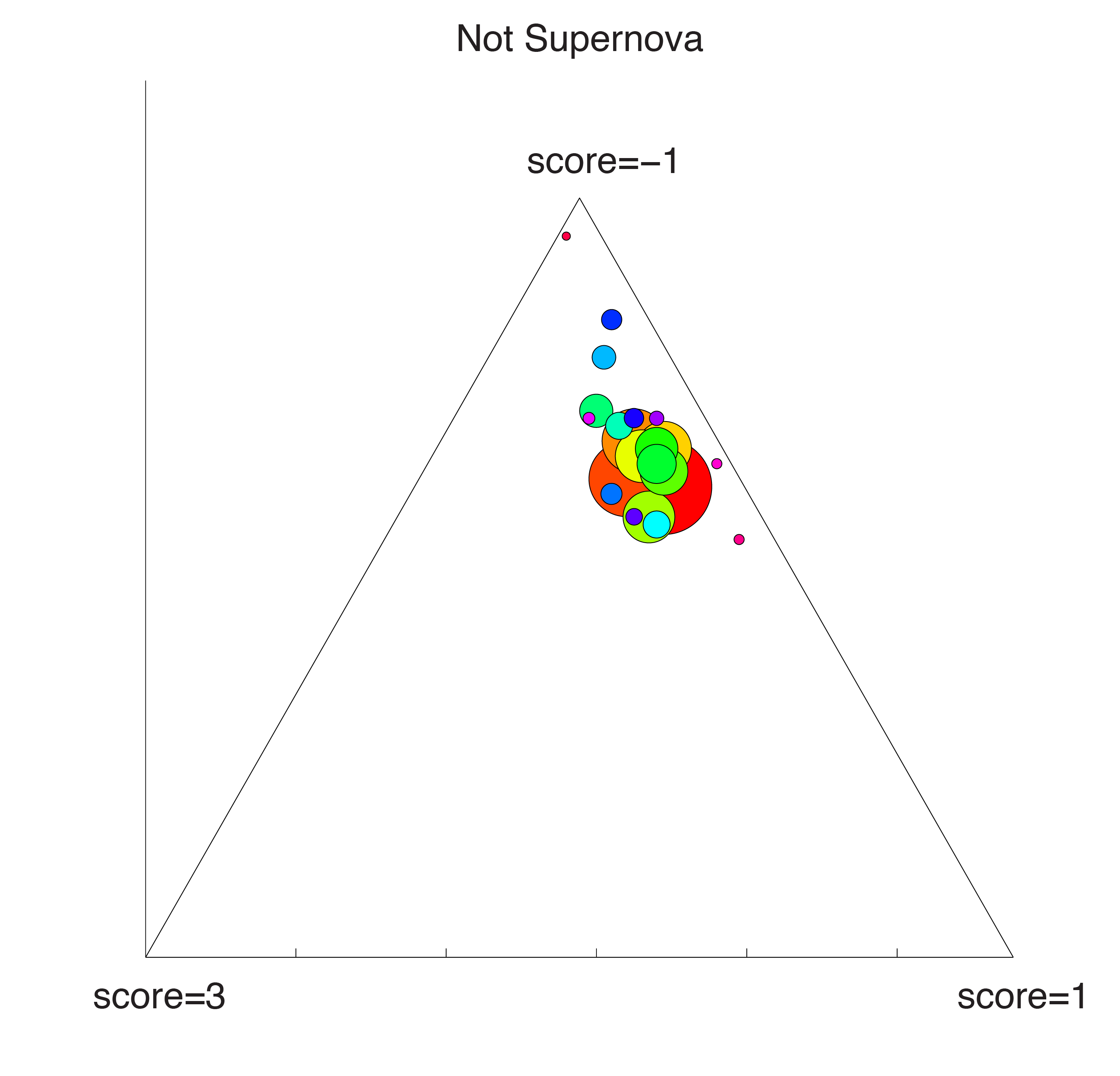}\includegraphics[clip=true,trim=0mm 0mm 0mm 0mm,width=0.4\textwidth,trim=15mm 0mm 0mm 0mm]{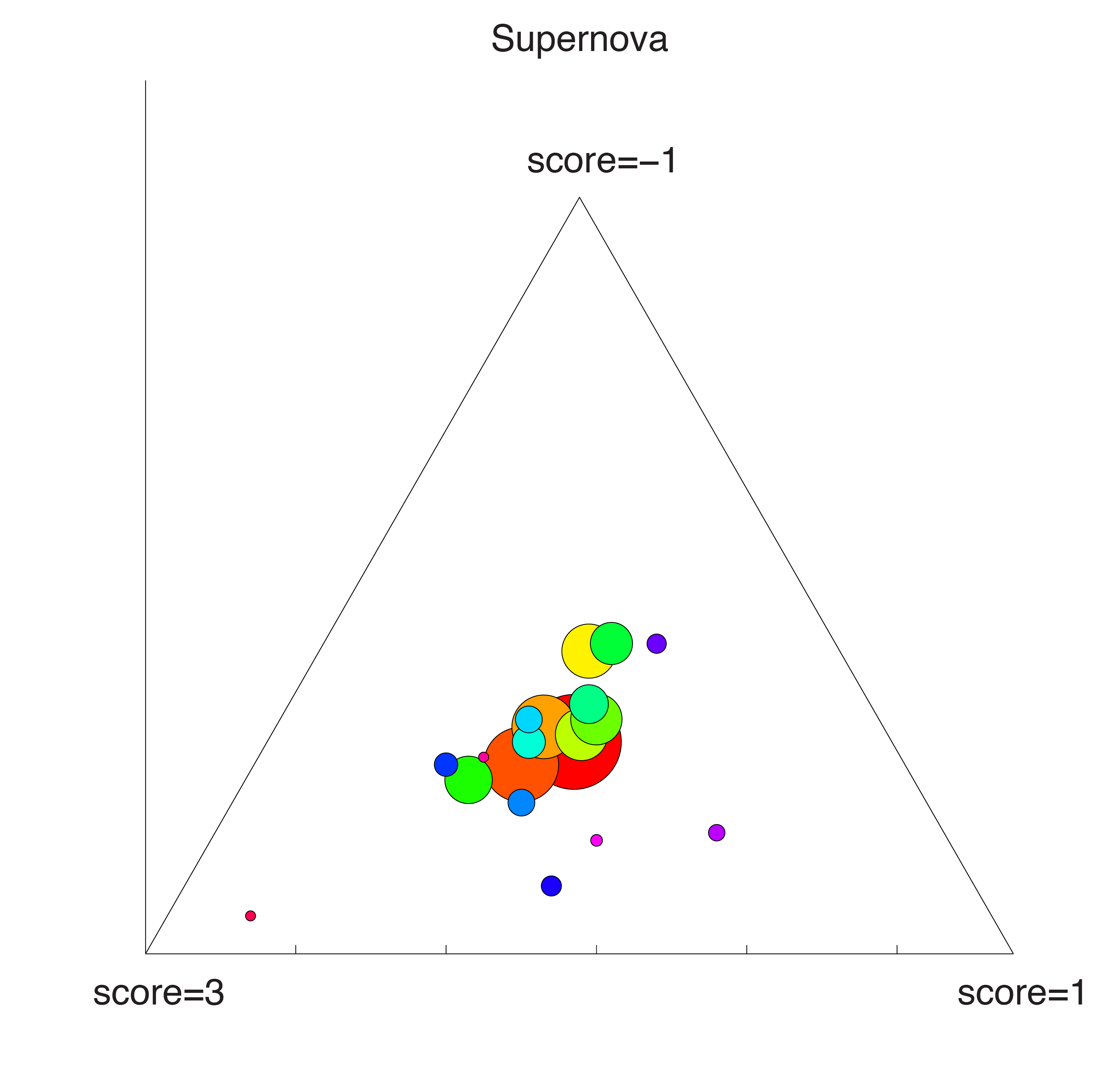}
\caption{Common task communities: ternary plot of means of expected confusion matrices for community members after 50,000 observations. Each point corresponds to one common task community and represents the mean of $\mathbb{E}[\boldsymbol{\pi}_j]$ for the community members. Proximity of a point to a vertex indicates the probability of outputting a particular score when presented with an object of true class ``supernova'' or ``not supernova''.}\label{fig:tasksCommunityDynamicsMeans200000}
\end{figure}

\begin{figure} [!ht]
\centering{} 
\includegraphics[clip=true,trim=0mm 0mm 0mm 0mm,width=0.4\textwidth,trim=15mm 0mm 0mm 0mm]{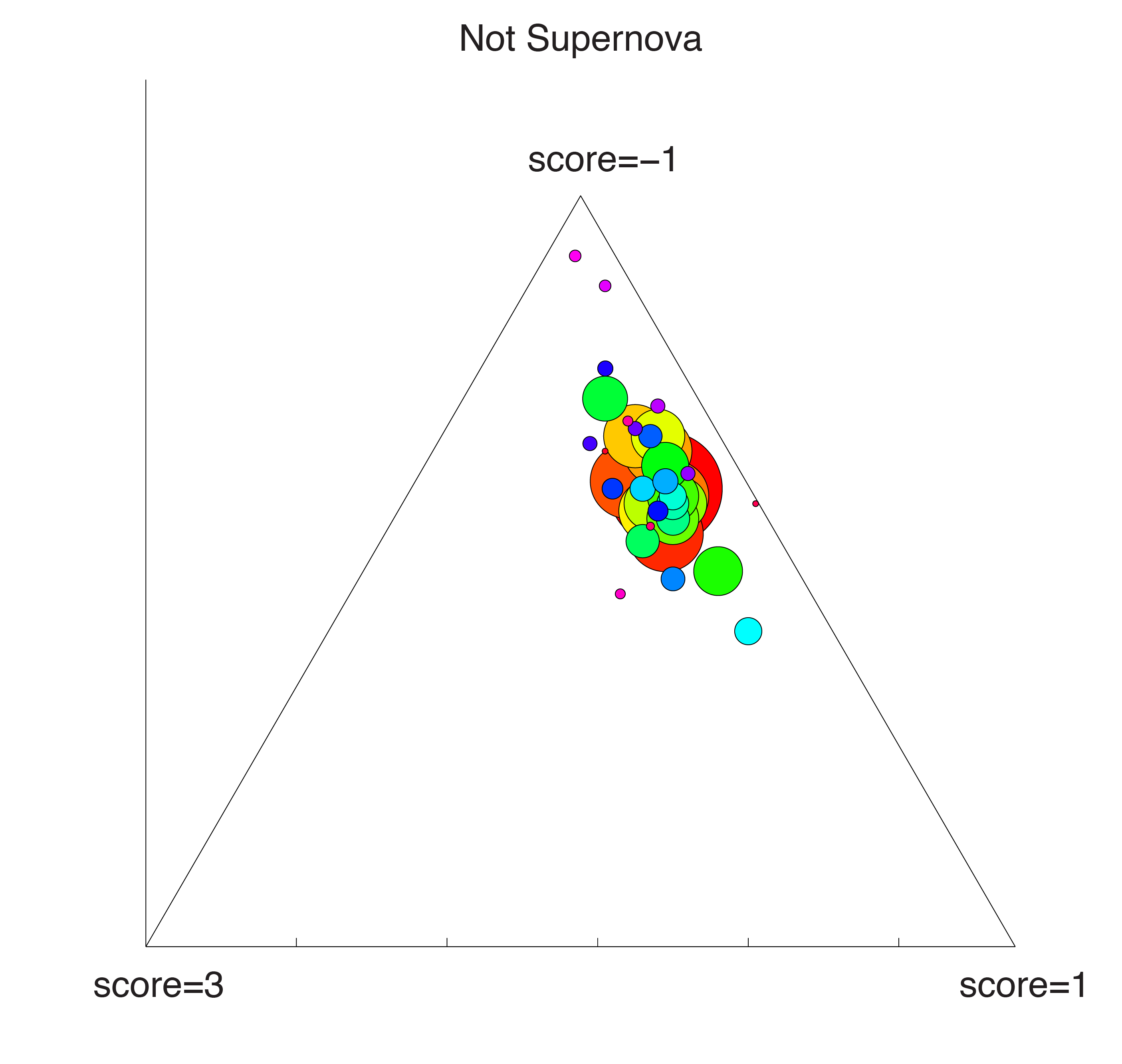}\includegraphics[clip=true,trim=0mm 0mm 0mm 0mm,width=0.4\textwidth,trim=15mm 0mm 0mm 0mm]{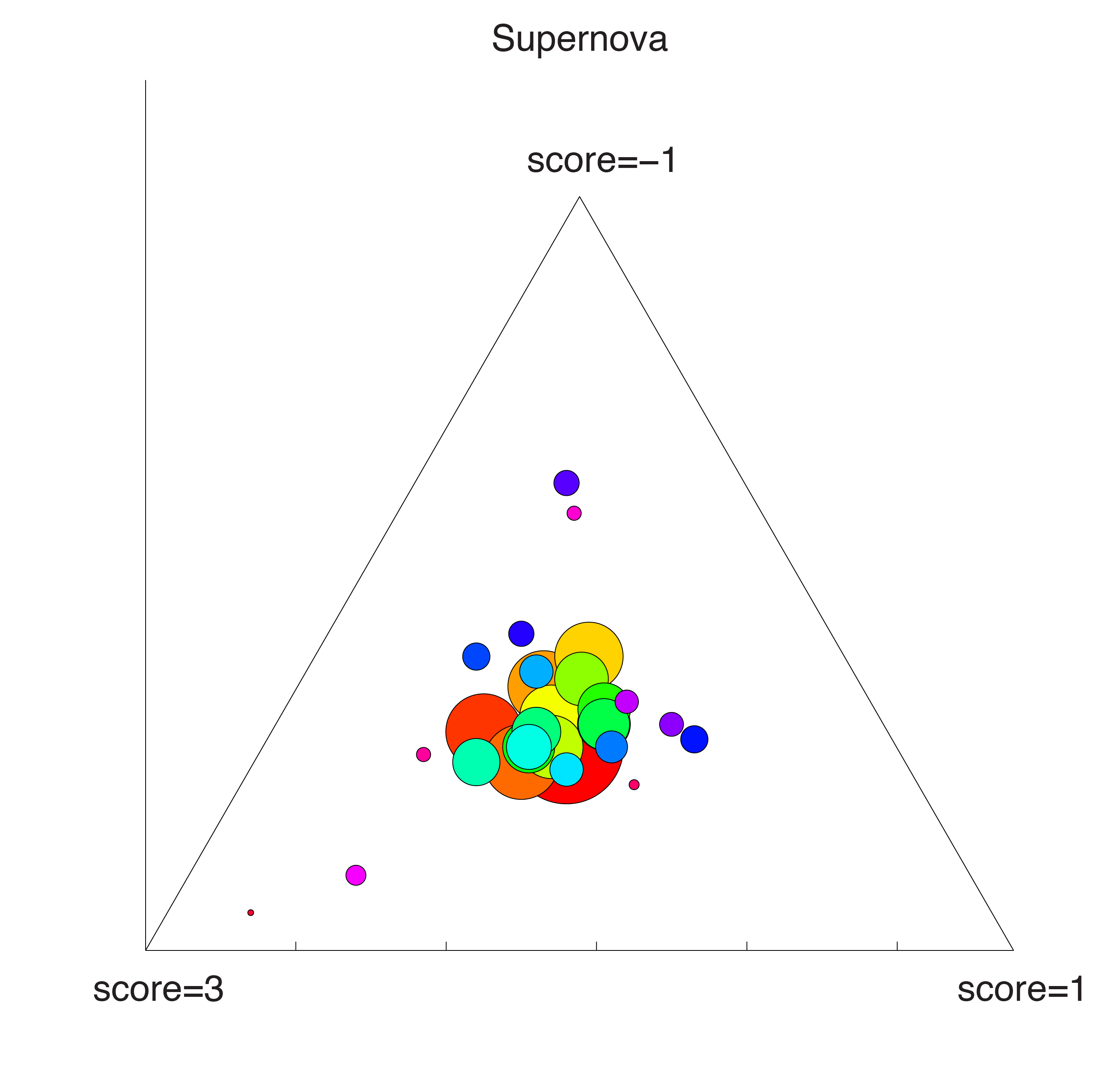}
\caption{Common task communities: ternary plot of means of expected confusion matrices for community members after 50,000 observations. Each point corresponds to one common task community and represents the mean of $\mathbb{E}[\boldsymbol{\pi}_j]$ for the community members. Proximity of a point to a vertex indicates the probability of outputting a particular score when presented with an object of true class ``supernova'' or ``not supernova''.}\label{fig:tasksCommunityDynamicsMeans493048}
\end{figure}

Finally, we look at the evolution of the common task communities to observe the effect of recent tasks on the community structure and confusion matrices. We wish to observe whether distinct communities are persistent as more tasks are completed. Changes to the structure inferred may be a result of observing more data about the base classifiers. Alternatively, individual behaviours may evolve as a result of learning new types of tasks or changes to individual circumstances, such as when a volunteer is available to carry out tasks. Our choice of community analysis method, given in Section \ref{sec:picomms} has the advantage that only a maximum number of communities need be chosen by the programmer, with the algorithm itself finding the most likely number of communities from the network itself.
Here we show the changes that occur in Galaxy Zoo Supernovae. We generated three co-occurrence networks from all tasks completed up to 50,000, 200,000 and 493,048 observations. As before, we remove base classifiers with fewer than 10 classifications to filter out edges that may constitute noise rather than significant similarity. The algorithm produced community structures with modularities of 0.67, 0.69 and 0.75 respectively, showing that good community structure is present for smaller periods of observations (see Section \ref{sec:picomms} for definition of modularity). Figures (\ref{fig:tasksCommunityDynamicsMeans50000}), (\ref{fig:tasksCommunityDynamicsMeans200000}) and (\ref{fig:tasksCommunityDynamicsMeans493048}) show the means of the community members at each of these time slices, weighted by node participation. Since DynIBCC models the dynamics of base classifier confusion matrices as a random walk, the observations closest to the current time have the 
strongest effect on the distribution over the confusion matrices. Therefore, the expected confusion matrices can readily be used to 
characterise a community at a given point in time. When calculating the means we use the expected confusion matrix from the most recent time-step for that network. 

In all three networks there is a persistent core for both true classes, where the means for the large communities remain similar. Some communities within this group move a small amount, for example, the large red community in the ``Supernova'' class. In contrast, we see more scattered small communities appear after 200,000 observations and at 493,048 observations. It is possible that the increase in number of base classifiers as we see more data means that previous individual outliers are now able to form communities with similar outliers. Therefore outlying communities could be hard to detect with smaller datasets. Many of these appear in the same place in only one of the figures, suggesting that they may contain new base classifiers that have made few classifications up to that point. Some are less transient however: the top-most community in the ``not supernova'' class in Figures (\ref{fig:tasksCommunityDynamicsMeans200000}) and (\ref{fig:tasksCommunityDynamicsMeans493048}) moves only a small amount. 
Similar sets of tasks may produce more extreme confusion matrices such as these for different agents at different times, implying that these tasks induce a particular bias in the confusion matrices. The changes we observe in Figures (\ref{fig:tasksCommunityDynamicsMeans50000}), (\ref{fig:tasksCommunityDynamicsMeans200000}) and (\ref{fig:tasksCommunityDynamicsMeans493048}) demonstrate how we can begin to identify the effect of different tasks on our view of the base classifiers by evaluating changes to the community structure after classifying certain objects. Future investigations may consider the need to modify the co-occurrence network to discount older task-based classifier associations. 

\section{Discussion} \label{sec:discussion}

In this paper we present a very computationally efficient, variational Bayesian, approach to imperfect multiple classifier combination. We evaluated the method using real data from the Galaxy Zoo Supernovae citizen science project, with 963 classification tasks, 1705 base classifiers and 26,558 observations. In our experiments, our method far outperformed all other methods, including weighted sum and weighted majority, both of which are often advocated as they also learn weightings for the base classifiers. For our variational Bayes method the required computational overheads were far lower than those of Gibbs sampling approaches, giving much shorter compute time, which is particularly important for applications that need to make regular updates as new data is observed, such as our application here. Furthermore, on this data set at least, the accuracy of predictions was also better than the slower sampling-based method. We have shown that social network analysis can be used to extract sensible structure from 
the pool of decision makers using information inferred by Bayesian classifier combination or task co-occurrence networks. This structure provides a useful grouping of individuals and gives valuable information about their decision-making behaviour. We extended our model to allow for on-line dynamic analysis and showed how this enables us to track the changes in time associated with individual base classifiers. We also demonstrated how the community structures change over time, showing the use of the dynamic model to update information about group members. 

Our current work considers how the rich information learned using our models can be exploited to improve the base classifiers, namely the human volunteer users. For example, we can use the confusion matrices, $\boldsymbol{\Pi}$, and the community structure to identify users who would benefit from more training. This could take place through interaction with user groups who perform more accurate decision making, for example via extensions of \emph{apprenticeship learning} \cite{abbeel_apprenticeship_2004}. We also consider ways of producing user specialisation via selective object presentation such that the overall performance of the human-agent collective is maximised. We note that this latter concept bears the hallmarks of \emph{computational mechanism design} \cite{dash_computational-mechanism_2003} and the incorporation of incentive engineering and coordination mechanisms into the model is one of our present challenges. Future work will also investigate selecting individuals for a task to maximise both 
our knowledge of the true labels and of the confusion matrices, for example, by looking at the effects of previous tasks on the confusion matrices. To bring these different aspects together, we consider a global utility function for a set of classification and training tasks indexed $i=1,...,N$ assigned to a set of base classifiers $k=1,...,K$. Classifiers assigned to object $i$ are part of coalition $\mathbf{C}_i$ to maximise the total expected value of these assignments:
\begin{equation}
 V(\mathbf{C}_1,...,\mathbf{C}_N) = \sum_{i=1}^N V_{object}(\mathbf{C}_i) + V_{dm}(\mathbf{C}_i) + V_{cost}(\mathbf{C}_i) \label{eq:weakcontrolvalues}
\end{equation}
where $V_{object}(\mathbf{C}_i)$ is the expected information gain about the true class of object $i$ from the classifiers in $\mathbf{C}_i$, $V_{dm}(\mathbf{C}_i)$ is the improvement to the decision makers through this assignment and $V_{cost}(k,i)$ captures other costs, such as payments to a decision maker. The value $V_{object}(\mathbf{C}_i)$ should be higher for classifiers in $\mathbf{C}_i$ that are independent, so coalitions of decision makers from different communities may be favoured as different experience and confusion matrices may indicate correlation is less likely. $V_{object}(\mathbf{C}_i)$ should also account for specialisations, for example, by members of the same common task community. $V_{dm}(\mathbf{C}_i)$, captures expected changes to confusion matrices that result from the combiner learning more about base classifiers and from base classifiers improving through training or experience. In Galaxy Zoo Supernovae, for example, the contributors are attempting to identify objects visually from 
a textual description. The description may leave some ambiguity, e.g. ``is the candidate roughly centred''. Seeing a range of images may alter how ``roughly'' the candidate can be centred before the contributor answers ``yes''. Thus the value $V_{dm}(k,i)$ will depend on the objects previously classified by classifier $k$. A key direction for future work is defining these values so that systems such as Galaxy Zoo Supernovae can feed back information from confusion matrices and community structure to improve the overall performance and efficiency of the pool of decision makers. Common task communities and $\boldsymbol{\pi}$ communities may play a central role in estimating the effects of task assignments and training on related individuals. The could also be exploited to reduce the size of the task assignment problem to one of choosing classifiers from a small number of groups rather than evaluating each classifier individually.

\section{Acknowledgements}

The authors would like to thank Chris Lintott at Zooniverse. We gratefully acknowledge funding from the UK Research Council EPSRC for project ORCHID, grant EP/I011587/1. Ioannis Psorakis is funded by Microsoft Research European PhD Scholarship Programme. 

\bibliographystyle{unsrt}
\bibliography{bcc_ext_revised}
\end{document}